\def \NN {\mathbb N}
\def \CC {\mathbb C}
\def \QQ {\mathbb Q}
\def \RR {\mathbb R}

\def \sameorder{\asymp}
\def \epsilon{\varepsilon}

\def \I  {{\mathcal I}}
\def \J  {{\mathcal J}}
\def \LL {{\mathcal L}}

\def \X  {{\mathcal X}}
\def \d {\text{d}}

\def \fine {{\hfill \qedsymbol}}

\renewcommand{\S}{{\mathcal S}}

\documentclass[12pt,reqno]{amsart}

\usepackage{amsmath,amssymb}

\numberwithin{equation}{section}

\begin{document}


\title[]{Twists, Euler products and a converse theorem for $L$-functions of degree 2 in the Selberg class}
\author[]{J.Kaczorowski  \lowercase{and} A.Perelli}
\date{}
\maketitle

\bigskip
{\bf Abstract.} We prove a general result relating the shape of the Euler product of an $L$-function to the analytic properties of certain linear twists of the $L$-function itself. Then, by a sharp form of the transformation formula for linear twists, we check the required analytic properties in the case of $L$-functions of degree 2 and conductor 1 in the Selberg class. Finally we prove a converse theorem, showing that $\zeta(s)^2$ is the only member of the Selberg class satisfying the above conditions and, moreover, having a pole at $s=1$.



\vskip.8cm
\section{Introduction}

\smallskip
In this paper we deal mainly with the $L$-functions of degree 2 and conductor 1 from the Selberg class $\S$. We refer to Selberg \cite{Sel/1989}, Conrey-Ghosh \cite{Co-Gh/1993}, to our survey papers \cite{Ka-Pe/1999b}, \cite{Kac/2006}, \cite{Per/2005}, \cite{Per/2004}, \cite{Per/2010} and to our forthcoming book \cite{book} for the basic information and results on the class $\S$ and on the extended Selberg class $\S^\sharp$ of $L$-functions. Moreover, we refer to the beginning of Section 3 for the definition of such classes, in particular for the data $\omega$, $Q$, $r$, $\lambda_j$ and $\mu_j$ which appear in the functional equation; see \eqref{3-0} below. Here we recall that degree and conductor of $F\in\S^\sharp$ are defined respectively by
\[
d_F = 2\sum_{j=1}^r\lambda_j \hskip2cm q_F= (2\pi)^{d_F}Q^2\prod_{j=1}^r\lambda_j^{2\lambda_j},
\]
that $\S_d$ and $\S^\sharp_d$ denote respectively the subclasses of $\S$ and $\S^\sharp$ of the functions of degree $d$ and that the Euler product of $F\in\S$ has the general form
\[
F(s) = \prod_p F_p(s),
\]
$F_p(s)$ being the $p$-th local factor of $F(s)$.

\medskip
The structure of the classes $\S$ and $\S^\sharp$ has been fully described for degrees $0\leq d<2$ in Conrey-Ghosh \cite{Co-Gh/1993} ($\S_d$ for $0\leq d <1$) and in our papers \cite{Ka-Pe/1999a} and \cite{Ka-Pe/2011} ($\S_1$ and $\S_d^\sharp$ for $0\leq d<2$). These results confirm the current conjectures on the structure of $\S$, i.e. the degree conjecture and the general converse problem. In particular, it turns out that $\S_d=\emptyset$ for $0<d<1$ and $1<d<2$, while $\S_1$ coincides with the GL$_1(\QQ)$ automorphic $L$-functions. In this paper we start investigating the next step, namely the description of $\S_2$. It is expected that $\S_2$ coincides with the GL$_2(\QQ)$ automorphic $L$-functions. Here we prove a rather special converse theorem, see Theorem 4 below, but the new ideas we employ appear to be suitable for further developments.

\medskip
Although the emphasis of the paper is mainly on degree 2 $L$-functions, we start with a general result relating linear twists and Euler product of the functions in $\S$. We recall that given $\alpha\in\RR$ and a Dirichlet character $\chi$, the linear (or additive) twist and the (multiplicative) twist of a function $F\in\S^\sharp$ with Dirichlet coefficients $a(n)$ are defined respectively by (as usual $e(x) = e^{2\pi i x}$)
\[
F(s,\alpha) = \sum_{n=1}^\infty\frac{a(n)e(-n\alpha)}{n^s} \hskip2cm F(s,\chi) =  \sum_{n=1}^\infty\frac{a(n)\chi(n)}{n^s}.
\]
Moreover, let
\[
N_F(\sigma,T) = |\{\rho=\beta+i\gamma: F(\rho)=0, \beta>\sigma, |\gamma|\leq T\}|
\]
be the zero-density function of $F(s)$. Further, for real numbers $d,h>0$ we define the class $M(d,h)$ as follows: $f\in M(d,h)$ if

\smallskip
\noindent
{\sl i)} $f(s)$ is meromorphic over $\CC$ and holomorphic for $\sigma<1$,

\noindent
{\sl ii)} for every $A<B$ there exists a constant $C=C(A,B)$ such that
\[
f(\sigma+it) \ll |\sigma|^{d|\sigma|}\big(\frac{h}{(2\pi e)^d}\big)^{|\sigma|} |\sigma|^C
\]
as $\sigma\to-\infty$ uniformly for $A\leq t\leq B$; the implied constant may depend on $d,h,A,B$ and $f(s)$. 

\smallskip
\noindent
Note that a function $F\in\S_d^\sharp$ belongs to $M(d,q_F)$, see Lemma 2.1 below. 

\bigskip
{\bf Theorem 1.} {\sl Let $d>0$ and $F\in\S_d$ be such that $N_F(\sigma,T)=o(T)$ for any fixed $\sigma>1/2$. Moreover, let $p$ be a prime number and $h>0$. Then the following statements are equivalent:

\smallskip
\noindent
i) for every $a$} (mod $p$), $a\not\equiv 0$ (mod $p$), $F(s,a/p)$ {\sl belongs to} $M(d,h)$;

\noindent
{\sl ii) for every $\chi$} (mod $p$), $\chi\neq \chi_0$, $F(s,\chi)$ {\sl belongs to $M(d,h)$ and}
\[
F_p(s) = \prod_{j=1}^{\partial_p} \big(1-\frac{\alpha_j(p)}{p^s}\big)^{-1} \quad \text{{\sl with $|\alpha_j(p)|\leq 1$ and} $\partial_p \leq \left[\frac{\log(h/q_F)}{\log p}\right]$.}
\]

\bigskip
{\bf Remarks.} {\bf 1.} Clearly, the interesting part of Theorem 1 is the implication {\sl i)} $\Rightarrow$ {\sl ii)}, since it shows that suitable analytic properties of the linear twists provide information on the shape of the Euler product. This phenomenon appears to be new and has no counterpart in the theory of classical $L$-functions, where the shape of the Euler product is given, in most cases, essentially by definition. Somehow, the $p$-th Euler factor is a measure of the difference between the groups of the additive and multiplicative characters (mod $p$). Indeed, the transition from multiplicative to additive twists is possible only when the $p$-th local factor is taken into account. Note, moreover, that the analytic properties required in Theorem 1 involve bounds on horizontal strips, while usually bounds on vertical strips are required in analytic number theory.

{\bf 2.} The form of the $p$-th local factor given by Theorem 1 is very close to the expected one. Indeed, it is expected that every $L$-function in $\S$ has local factors of the form
\[
F_p(s) = \prod_{j=1}^{\partial_F}\big(1-\frac{\alpha_j(p)}{p^s}\big)^{-1}
\]
with $|\alpha_j(p)|\leq 1$ for all $j$ and $p$, and $|\alpha_j(p)|=1$ for $j=1,\dots,\partial_F$ and all but finitely many primes $p$. Moreover, it is expected that $\partial_F = d_F$; see \cite{Ka-Pe/2008} for a discussion of these topics. On the other hand, it is also expected that $q_F\in\NN$ for every $F\in\S$, and that for any primitive $\chi$ (mod $m$) with $(m,q_F)=1$ the twist $F(s,\chi)$ belongs to $\S$ and has conductor $q_{F^\chi}=q_Fm^{d_F}$; see \cite{Ka-Pe/2000} for a discussion of these topics. As a consequence, we expect that for every $p$ the above linear twists $F(s,a/p)$ belong to the class $M(d_F, q_Fp^{d_F})$, in which case Theorem 1 implies that $\partial_p \leq [d_F]$ for all $p$.

{\bf 3.} The hypothesis on $N_F(\sigma,T)$ is a weak zero-density bound, which is known to hold for classical $L$-functions of degree 1 and 2; see the Corollary below for a general version in $\S_2$. It is a challenging problem to extend such a bound to every degree, even for classical $L$-functions. Note that, although the bound is just a little better than trivial, it is nevertheless quite interesting, since possibly it characterizes the $L$-functions satisfying the Riemann Hypothesis inside a rather large class of Dirichlet series with functional equation; see Kaczorowski-Kulas \cite{Ka-Ku/2007}.

\bigskip
From now on we shall deal with degree 2 $L$-functions. Our first aim is to show that statement {\sl i)} of Theorem 1 holds for $L$-functions of degree 2 and conductor 1; see Theorem 3 below. The main tool in such a proof is a precise form of the transformation formula for linear twists obtained in Lemma A of \cite{Ka-Pe/2002b}, see also \cite{Ka-Pe/2010}. Writing 
\[
2\sum_{j=1}^r(\mu_j-1/2) = \eta_F+i\theta_F \hskip2cm \eta_F,\theta_F\in\RR,
\]
\[
 \omega_F^* = \omega e^{-i\frac{\pi}{2}(\eta_F+1)}\big(\frac{q_F}{(2\pi)^{2}}\big)^{i\frac{\theta_F}{2}}\prod_{j=1}^r\lambda_j^{-2i\Im \mu_j}, 
\hskip1cm \bar{F}(s,\alpha)= \sum_{n=1}^\infty \frac{\overline{a(n)}e(-n\alpha)}{n^s}
\]
respectively for the weight, the shift, the root number and the linear twist of the conjugate of $F\in\S^\sharp_2$, we have

\bigskip
{\bf Theorem 2.} {\sl Let $F\in\S_2^\sharp$ and $\alpha>0$. Then for every integer $K>0$ there exist polynomials $Q_0(s),...,Q_K(s)$, with $Q_0(s)\equiv 1$, such that
\[
F(s,\alpha) = -i\omega_F^* (\sqrt{q_F}\alpha)^{2s-1+i\theta_F} \sum_{\nu=0}^K \big(\frac{iq_F\alpha}{2\pi}\big)^\nu Q_\nu(s) \bar{F}(s+\nu+i\theta_F,-\frac{1}{q_F\alpha}) + H_K(s,\alpha).
\]
Here $H_K(s,\alpha)$ is holomorphic for $-K+1/2<\sigma<2$ and $|s|<2K$, and satisfies
\[
H_K(s,\alpha) \ll (AK)^K
\]
with a suitable constant $A=A(F,\alpha)>0$. Moreover, $\deg Q_\nu(s)=2\nu$ and}
\[
\begin{split}
Q_\nu(s) &\ll \frac{(A(|s|+1))^{2\nu}}{\nu!} \hskip2cm 1\leq \nu\leq \min(|s|,K) \\
Q_\nu(s) &\ll (AK)^K \hskip3.2cm |s|\leq 2K, \, \nu\leq K.
\end{split}
\]

\bigskip
The proof of Theorem 2 is rather complicated and is given in Section 3. Theorem 2 is the key ingredient in the proof of 

\bigskip
{\bf Theorem 3.} {\sl Let $F\in\S^\sharp_2$ with $q_F=1$. Then for every $q\geq 1$ and $1\leq a\leq q$ with $(a,q)=1$ the linear twist $F(s,a/q)$ belongs to $M(2,q^2)$.}

\bigskip
From Theorems 1 and 3 we get

\bigskip
{\bf Corollary.} {\sl Let $F\in\S_2$ with $q_F=1$. Then for every prime $p$ and every} $\chi$ (mod $p$), $\chi\neq \chi_0$, {\sl the twist $F(s,\chi)$ belongs to $M(2,p^2)$ and}
\[
F_p(s) = \prod_{j=1}^{\partial_p} \big(1-\frac{\alpha_j(p)}{p^s}\big)^{-1} \quad \text{{\sl with $|\alpha_j(p)|\leq 1$ and} $\partial_p \leq 2$.}
\]

\bigskip
Note that the same arguments used in Theorem 3 and its Corollary give that $(s-1)^{m_F}F(s,\chi)$ is entire for every $\chi$ (mod $p$), $\chi\neq \chi_0$, where $m_F$ is the order of pole of $F(s)$ at $s=1$.

\bigskip
The above Corollary allows us to get a converse theorem for $\S_2$. In his famous paper \cite{Maa/1949}, Maass showed, among other, that the vector space of the functions $F\in\S^\sharp$ satisfying the functional equation of $\zeta(s)^2$ is 1-dimensional, and hence generated by $\zeta(s)^2$. From the Corollary we deduce the following characterization of $\zeta(s)^2$.

\bigskip
{\bf Theorem 4.} {\sl Let $F\in\S_2$ with $q_F=1$ and a pole at $s=1$. Then $F(s) = \zeta(s)^2$}.

\bigskip
Comparing Maass' converse theorem with Theorem 4, we see that the main difference is that we deal with a general degree 2 functional equation but we assume that $F(s)$ has an Euler product, while Maass does not need the Euler product but deals with a special degree 2 functional equation.

\medskip
{\bf Remark.} We finally note that our converse theorem is proved by showing that under the hypotheses of Theorem 4, the Euler product of $F(s)$ coincides with the Euler product of $\zeta(s)^2$. Twists are needed to prove the required properties of the Euler product and, as we already outlined after Theorem 1, are used in a definitely different way with respect to the classical converse theorems of Weil's type (see Ch.7 of Iwaniec's book \cite{Iwa/1997}). Therefore, in some sense our result realizes another instance of the approach to converse theorems via Euler products proposed in the paper by Conrey-Farmer \cite{Co-Fa/1995} (see also Conrey-Farmer-Odgers-Snaith \cite{CFOS/2007}). Such an approach represents an interesting alternative to the classical converse theorems based on twists, and on Rankin-Selberg convolutions for higher degrees; see Cogdell and Piatetski-Shapiro \cite{Co-PS/2001}.

\bigskip
{\bf Acknowledgements.} We thank Sandro Bettin and Brian Conrey for carefully reading a previous version of this paper and suggesting several improvements in the presentation. This research was partially supported by the Istituto Nazionale di Alta Matematica, by a grant PRIN2008 and by grant N N201 605940 of the National Science Centre.

\bigskip
\section{Proof of Theorem 1}

\smallskip
Given $F\in\S^\sharp$ satisfying the functional equation in \eqref{3-0} below we define
\[
\tau_F = \max_{1\leq j\leq r}\big|\frac{\Im\mu_j}{\lambda_j}\big|.
\]
It is easy to check by means of the criteria in \cite{Ka-Pe/2000} that $\tau_F$ is an invariant. Recalling the definition of the class $M(d,h)$ in the Introduction, we have

\medskip
{\bf Lemma 2.1.} {\sl Let $F\in\S^\sharp_d$ with $d>0$; then $F\in M(d,q_F)$. Moreover, if $[A,B]\cap [-\tau_F,\tau_F] = \emptyset$ and $\sigma\geq1$ we also have
\begin{equation}
\label{2-0}
F(-\sigma+it) \gg \sigma^{d\sigma} \big(\frac{q_F}{(2\pi e)^d}\big)^\sigma \sigma^C
\end{equation}
for some $C=C(A,B)$, uniformly for $A\leq t\leq B$ as $\sigma\to+\infty$.}

\medskip
{\sl Proof.} In this proof we write $s=-\sigma+it$ and assume that $\sigma\geq 1$. The regularity conditions required by the class $M(d,q_F)$ are clearly satisfied. By the functional equation and the reflection formula for the $\Gamma$ function we obtain
\[
F(s) = \omega Q^{1-2s} S(s)G(s) \bar{F}(1-s)
\]
with
\[
S(s) = \pi^{-r} \prod_{j=1}^r \sin(\pi(\lambda_j s+\mu_j))      \hskip1cm
G(s) = \prod_{j=1}^r \Gamma(\lambda_j(1-s)+\mu_j) \Gamma(1-\lambda_js-\mu_j).
\]
Uniformly for $A\leq t\leq B$ we have $S(s) \ll 1$, and if $[A,B]\cap [-\tau_F,\tau_F] = \emptyset$ we also have $|S(s)| \sameorder 1$, where $f\sameorder g$ means $g \ll f \ll g$. Moreover, the bound $|\bar{F}(1-s)| \sameorder 1$ holds since $\sigma\geq 1$. Hence uniformly for $A\leq t\leq B$ we have
\begin{equation}
\label{2-1}
|F(s)| \ll Q^{2\sigma} |G(s)|,
\end{equation}
and if $[A,B]\cap [-\tau_F,\tau_F] = \emptyset$ we also have
\begin{equation}
\label{2-2}
|F(s)| \sameorder Q^{2\sigma} |G(s)|.
\end{equation}
Writing $\beta=\prod_{j=1}^r\lambda_j^{2\lambda_j}$, from Stirling's formula we get
\[
\log|G(s)| = \sum_{j=1}^r\big\{(2\lambda_j\sigma+O(1))\log(\lambda_j\sigma)-2\lambda_j\sigma +O(1)\big\} = d\sigma\log\sigma + (\log\beta -d)\sigma +O(\log\sigma),
\]
and hence
\begin{equation}
\label{2-3}
|G(s)| \sameorder \sigma^{d\sigma} \big(\frac{q_F}{(2\pi e)^d}\big)^\sigma Q^{-2\sigma} \sigma^{O(1)}.
\end{equation}
The lemma follows then from \eqref{2-1}-\eqref{2-3}. \fine

\medskip
{\bf Lemma 2.2.} {\sl Let $F\in\S_d$ with $d>0$. Let $p$ be a prime number, $\sigma>1$ and $\tau(\chi)$ denote the Gauss sum. For any $\chi$} (mod $p$), $\chi\neq \chi_0$, {\sl we have
\begin{equation}
\label{2-4}
F(s,\chi) = \frac{1}{\tau(\bar{\chi})} \sum_{a=1}^p \bar{\chi}(a) F(s,-a/p),
\end{equation}
while for any $(a,p)=1$ we have}
\begin{equation}
\label{2-5}
F(s,-a/p) = \frac{1}{p-1} \sum_{\substack{\chi \, \text{(mod $p$)} \\ \chi\neq \chi_0}} \chi(a)\tau(\bar{\chi}) F(s,\chi) - \big(\frac{p}{p-1}\frac{1}{F_p(s)}-1\big)F(s).
\end{equation}

\medskip
{\sl Proof.} Equation \eqref{2-4} is standard since $\chi$ is primitive. From \eqref{2-4}, writing $\delta_{a,b}=1$ if $a\equiv b$ (mod $p$) and $\delta_{a,b}=0$ otherwise, we get
\[
\begin{split}
\frac{1}{p-1} \sum_{\substack{\chi \, \text{(mod $p$)} \\ \chi\neq \chi_0}} \chi(a)\tau(\bar{\chi}) &F(s,\chi) = \frac{1}{p-1} \sum_{\substack{\chi \, \text{(mod $p$)} \\ \chi\neq \chi_0}} \sum_{b=1}^p \chi(a) \bar{\chi}(b) F(s,-b/p) \\
&= \sum_{b=1}^{p-1}\big(\frac{1}{p-1} \sum_{\chi \, \text{(mod $p$)}} \chi(a) \bar{\chi}(b) - \frac{1}{p-1}\big) F(s,-b/p) \\
&= \sum_{b=1}^{p-1}\big(\delta_{a,b} - \frac{1}{p-1}\big) F(s,-b/p) = F(s,-a/p) - \frac{1}{p-1} \sum_{b=1}^{p-1} F(s,-b/p) \\
&= F(s,-a/p) - \frac{1}{p-1} \sum_{n=1}^\infty \frac{a(n)}{n^s}\sum_{b=1}^{p-1} e(\frac{bn}{p}) \\
&= F(s,-a/p) - \sum_{p|n}\frac{a(n)}{n^s} - \frac{1}{p-1} \sum_{p\nmid n} \frac{a(n)}{n^s} \big(\sum_{b=0}^{p-1} e(\frac{bn}{p}) -1\big) \\
&= F(s,-a/p) -F(s) + \frac{F(s)}{F_p(s)} +\frac{1}{p-1} \frac{F(s)}{F_p(s)},
\end{split}
\]
and the lemma follows. \fine

\medskip
Now we are ready for the proof of Theorem 1. We first note that $F(s,a/p)$ belongs to $M(d,h)$ for every $a$ (mod $p$), $a\not\equiv 0$ (mod $p$), if and only if $F(s,-a/p)$ satisfies the same conditions, since $F(s,-a/p)=F(s,(p-a)/p)$. We start with the proof that {\sl ii)} $\Rightarrow$ {\sl i)}. From \eqref{2-5} in Lemma 2.2 and our assumptions we have that $F(s,a/p)$ is meromorphic on $\CC$, holomorphic for $\sigma<1$ and satisfies
\[
F(s,-a/p) \ll \max_{\chi\neq \chi_0}|F(s,\chi)| + \big|\frac{F(s)}{F_p(s)}\big| + |F(s)| \ll \max_{\chi\neq \chi_0}|F(s,\chi)| + |F(s)|p^{|\sigma|\partial_p}.
\]
Therefore by Lemma 2.1 we have, as $\sigma\to-\infty$ uniformly for $A\leq t\leq B$,
\[
\begin{split}
F(s,-a/p) &\ll |\sigma|^{d|\sigma|} \big(\frac{h}{(2\pi e)^d}\big)^{|\sigma|} |\sigma|^C + |\sigma|^{d|\sigma|} \big(\frac{q_F}{(2\pi e)^d}\big)^{|\sigma|} \big(\frac{h}{q_F}\big)^{|\sigma|} |\sigma|^C \\
&\ll |\sigma|^{d|\sigma|} \big(\frac{h}{(2\pi e)^d}\big)^{|\sigma|} |\sigma|^C
\end{split}
\]
with some constant $C$. Hence $F(s,a/p)$ belongs to $M(d,h)$, and the implication is proved.

\medskip
To prove that {\sl i)} $\Rightarrow$ {\sl ii)} we first note that by \eqref{2-4} in Lemma 2.2 and our assumptions we have, for $\chi\neq \chi_0$, that $F(s,\chi)$ is meromorphic on $\CC$, holomorphic for $\sigma<1$ and satisfies
\[
F(\sigma+it,\chi) \ll \max_{a \, \text{(mod $p$)}} |F(\sigma+it,-a/p)| \ll |\sigma|^{d|\sigma|} \big(\frac{h}{(2\pi e)^d}\big)^{|\sigma|} |\sigma|^C
\]
with some constant $C$, as $\sigma\to-\infty$ uniformly for $A\leq t\leq B$. Hence $F(s,\chi)$, $\chi\neq \chi_0$, belongs to $M(d,h)$. Concerning the $p$-th Euler factor $F_p(s)$, it is easy to see that $1/F_p(s)$ is holomorphic for $\sigma>\vartheta$ (see our survey \cite{Ka-Pe/1999b} and the beginning of Sect.3 below for the definition of $\vartheta$), and by \eqref{2-5} in Lemma 2.2 it is meromorphic on $\CC$. Moreover, since $F(s)\neq 0$ in the half-plane $\sigma<0$ apart from the trivial zeros, thanks to the $\frac{2\pi i}{\log p}$-periodicity we have that $1/F_p(s)$ is holomorphic in the same half-plane. As a consequence, the singularities of $1/F_p(s)$ may come only from the zeros of $F(s)$ in the strip $0\leq \sigma\leq \vartheta$. But our hypothesis on $N_F(\sigma,T)$, the functional equation and the fact that $\vartheta<1/2$ imply that the number of such singularities up to $T$ is $o(T)$, hence in view of the $\frac{2\pi i}{\log p}$-periodicity we deduce that $1/F_p(s)$ is entire. In addition, we may write
\[
\frac{1}{F_p(s)} = E(p^{-s})
\]
with $E(z)$ entire. Let now $T_0>\tau_F$ and consider the strip $S$ with $\sigma\leq -1$ and $T_0\leq t\leq T_0 + \frac{2\pi}{\log p}$. From \eqref{2-5} of Lemma 2.2, \eqref{2-0} of Lemma 2.1 and the fact the twists $F(s,\chi)$ and $F(s,a/p)$ belong to $M(d,h)$ we have
\[
E(p^{-s}) = \frac{1}{F_p(s)} \ll \frac{|\sigma|^{d|\sigma|} \big(\frac{h}{(2\pi e)^d}\big)^{|\sigma|} |\sigma|^C}{|\sigma|^{d|\sigma|} \big(\frac{q_F}{(2\pi e)^d}\big)^{|\sigma|} |\sigma|^{C'}} \ll \big(\frac{h}{q_F}\big)^{|\sigma|} |\sigma|^{C''} \ll p^{|\sigma| (\frac{\log(h/q_F)}{\log p} + \epsilon)}
\]
for every $\epsilon>0$. Therefore $E(z)$ is a polynomial of degree $\leq [\frac{\log(h/q_F)}{\log p}]$, hence
\begin{equation}
\label{2-6}
F_p(s) = \prod_{j=1}^{\partial_p} \big(1-\frac{\alpha_j(p)}{p^s}\big)^{-1}
\end{equation}
with $\partial_p \leq [\frac{\log(h/q_F)}{\log p}]$. The last assertion, namely $|\alpha_j(p)|\leq 1$, is a consequence of the Ramanujan condition $a(n)\ll n^\epsilon$ for every $\epsilon>0$, where $a(n)$ are the coefficients of $F(s)$ (see the beginning of Sect.3 below). Indeed, writing
\[
G_p(z) = \sum_{m=0}^\infty a(p^m)z^m
\]
the Ramanujan condition implies that $G_p(z)$ is holomorphic for $|z|<1$. On the other hand, since $F_p(s) = G_p(p^{-s})$ by \eqref{2-6} we also have
\[
G_p(z) = \prod_{j=1}^{\partial_p}(1-\alpha_j(p)z)^{-1},
\]
with poles at the points $z=\alpha_j(p)^{-1}$. Thus $|\alpha_j(p)|\leq 1$ and Theorem 1 is proved. \fine

\medskip
{\bf Remark.} We wish to thank Giuseppe Molteni for the above elegant proof that $|\alpha_j(p)|\leq 1$.

\bigskip
\section{Proof of Theorem 2}

\smallskip
{\bf 1. Definitions and notation.} We start with the definition of $\S$ and $\S^\sharp$. $F\in\S$ if

\smallskip
\noindent
{\sl i)} $F(s)$ is an absolutely convergent Dirichlet series for $\sigma>1$; 

\noindent
{\sl ii)} $(s-1)^mF(s)$ is
an entire function of finite order for some integer $m\geq 0$; 

\noindent
{\sl iii)} $F(s)$ satisfies a functional equation of type $\Phi(s) = \omega\bar{\Phi}(1-s)$,
where $|\omega|=1$ and
\begin{equation}
\label{3-0}
\Phi(s) = Q^s\prod_{j=1}^r\Gamma(\lambda_js+\mu_j)F(s)
\end{equation}
with $r\geq 0$, $Q>0$, $\lambda_j>0$, $\Re\mu_j\geq 0$ (here and in the sequel we write
$\bar{f}(s) =\overline{f(\overline{s})}$);

\noindent
{\sl iv)} the Dirichlet coefficients $a(n)$ of $F(s)$ satisfy $a(n) \ll n^\epsilon$ for every
$\epsilon>0$; 

\noindent
{\sl v)} $\log F(s)$ is a Dirichlet series with coefficients $b(n)$ satisfying $b(n)=0$ unless
$n=p^m$, $m\geq 1$, and $b(n)\ll n^\vartheta$ for some $\vartheta<1/2$. 

\smallskip
\noindent
The extended Selberg class $\S^\sharp$ consists of the non-zero functions satisfying only axioms {\sl i)}, {\sl ii)} and {\sl iii)}.

\medskip
Let $F\in\S_2^\sharp$ and consider the $H$-invariants
\[
H_F(n) = 2\sum_{j=1}^r\frac{B_n(\mu_j)}{\lambda_j^{n-1}} \hskip2cm n=0,1,...
\]
where $B_n(x)$ is the $n$-th Bernoulli polynomial; see \cite{Ka-Pe/2002a} for properties of such invariants. Note that $H_F(0) = d_F$ is the degree and 
\[
H_F(1) = 2\sum_{j=1}^r (\mu_j-\frac12) = \xi_F = \eta_F + i\theta_F
\]
is the $\xi$-invariant (in the Introduction we already defined its real and imaginary parts). For $\nu,\mu=1,2,...$ we define the polynomials $R_\nu(s) = R_{\nu,F}(s)$ and $V_\mu(s)=V_{\mu,F}(s)$ as
\begin{equation}
\label{3-0bis}
\begin{split}
R_\nu(s) &= B_{\nu+1}(1-2s-i\theta_F) + B_{\nu+1}(1) \\
&+ \frac12  \sum_{k=0}^{\nu+1}{\nu+1\choose k}\big((-1)^\nu H_F(k)s^{\nu+1-k} - \overline{H_F(k)}(1-s)^{\nu+1-k}\big)
\end{split}
\end{equation}
\begin{equation}
\label{3-1}
V_\mu(s) = (-1)^\mu\sum_{m=1}^\mu \frac{1}{m!} \sum_{\substack{\nu_1\geq 1,...,\nu_m\geq 1\\ \nu_1+...+\nu_m=\mu}} \prod_{j=1}^m\frac{R_{\nu_j}(s)}{\nu_j(\nu_j+1)}.
\end{equation}
We also define $Q_0(s)\equiv1$ and, for $\nu=1,2,...$, the functions $Q_\nu(s) = Q_{\nu,F}(s)$ by means of the formula
\[
\exp\big(\sum_{\nu=1}^\infty \frac{(-1)^\nu R_\nu(s)}{\nu(\nu+1)} \frac{1}{(w+2s-1+i\theta_F)^\nu}\big) \approx 1 + \sum_{\nu=1}^\infty \frac{Q_\nu(s)}{(w-1)\cdots(w-\nu)}.
\]
Here $\approx$ means asymptotic expansion as $w\to\infty$, and the $Q_\nu(s)$'s turn out to be polynomials; see Lemma 3.15 below for more details. 

\medskip
We write $w=u+iv$ and, for a given $s$, define the contour $\LL(s)$ as follows:
\[
\LL(s) = \LL_{-\infty}(s) \cup \LL_{\infty}(s)
\]
where
\[
\LL_{-\infty}(s) = (-\sigma+c_0-i\infty,-\sigma+c_0+it_0] \cup [-\sigma+c_0+it_0,-\sigma-c_0+it_0]
\]
\[
\LL_\infty(s) = [-\sigma-c_0+it_0,-\sigma-c_0+i\infty).
\]
Here $t_0=t_0(s) = c_1|s|+c_2$, $c_0,c_1,c_2>0$ being sufficiently large constants depending on $F(s)$ to be chosen later on. Moreover, we denote by $\LL^*_{-\infty}(s)$ the half-line $1-2s-i\theta_F - \LL_\infty(s)$ taken with the positive orientation, hence
\[
\LL_{-\infty}^*(s) = (1-\sigma+c_0-i\infty, 1-\sigma+c_0-it_0^*]
\]
with $t_0^*=t_0^*(s) =t_0 +2t+\theta_F$. Further, we let
\[
\LL_\infty^*(s)=[1-\sigma+c_0-it_0^*, N+1] \cup [N+1,N+1+i\infty),
\]
where the positive integer $N$ will be chosen later on (see \eqref{P6bis} below), and write
\[
\LL^*(s) = \LL^*_{-\infty}(s) \cup \LL_\infty^*(s).
\]

\medskip
We let
\[
G(s,w) = \frac{(2\pi)^{1-r}}{\Gamma(1-w)} \prod_{j=1}^r \Gamma(\lambda_j(1-s-w)+\bar{\mu}_j) \Gamma(1-\lambda_j(s+w)-\mu_j),
\]
\[
S(s,w) = \frac{2^{r-1}}{\sin\pi w} \prod_{j=1}^r \sin\big(\pi(\lambda_j(s+w)+\mu_j\big).
\]
Moreover, we define the coefficients $C_{\mu,\ell}$, $\ell\geq \mu\geq 1$, and $A_{\mu,\nu}(s)$, $\nu\geq \mu\geq 1$, by the asymptotic expansions
\begin{equation}
\label{3-2}
\frac{1}{w^\mu} \approx \sum_{\ell=\mu}^\infty \frac{C_{\mu,\ell}}{(w-1)\cdots(w-\ell)}
\end{equation}
\begin{equation}
\label{3-3}
\frac{1}{(w+2s-1+i\theta_F)^\mu} \approx  \sum_{\nu=\mu}^\infty \frac{A_{\mu,\nu}(s)}{(w-1)\cdots(w-\nu)}.
\end{equation}
We refer to Lemmas 3.13 and 3.14 below for the precise meaning of \eqref{3-2} and \eqref{3-3}.

\medskip
Finally, $A,c,c',\tilde{c},c^*,...$ will denote positive constants, possibly depending on $F(s)$ (also via a dependence on $c_0,c_1,c_2$ above), not necessarily the same at each occurrence. The constants in the $\ll$- and $O$-symbols may also depend on $F(s)$ (again, also via $c_0,c_1,c_2$).

\bigskip
{\bf 2. Lemmas.} In the following lemmas we always assume that the parameters $\lambda_j,\mu_j,\theta_F$,... come from a function $F\in\S_2^\sharp$, unless otherwise specified.

\bigskip
{\bf Lemma 3.1.} {\sl For $w\in\LL_\infty(s)$ and $s\in\CC$ we have
\[
G(s,w) \ll e^{-\frac{\pi}{2}v} A^{|s|} (v+|\sigma|)^{|\sigma|+c} \hskip2cm v\to +\infty
\]
with suitable constants} $A,c>0$.

\bigskip
{\sl Proof.} For $w\in\LL_\infty(s)$ we have
\[
\Re(\lambda_j(1-s-w)+\bar{\mu}_j) = \lambda_j(1+c_0)+\Re\mu_j>0
\]
\[
\Im(\lambda_j(1-s-w)+\bar{\mu}_j) = -\lambda_j(t+v) - \Im\mu_j,
\]
hence by Stirling's formula
\begin{equation}
\label{L3-1/1}
\Gamma(\lambda_j(1-s-w)+\bar{\mu}_j) \ll e^{-\lambda_j\frac{\pi}{2}|t+v|} |t+v|^c \ll e^{-\lambda_j\frac{\pi}{2}v} A^{|s|} v^c.
\end{equation}
Similarly we have
\begin{equation}
\label{L3-1/2}
\Gamma(1-\lambda_j(s+w)-\mu_j) \ll e^{-\lambda_j\frac{\pi}{2}v} A^{|s|} v^c.
\end{equation}
In order to treat $1/\Gamma(1-w)$ we first note that $1-w$ belongs to a region where the Stirling formula is applicable, and we have
\[
\log|\Gamma(1-w)| = (\frac12-u)\log|1-w| +v\arg(1-w) +u+O(1).
\]
Moreover
\[
\arg(1-w) = -\frac{\pi}{2} + O(\arctan\frac{|1+\sigma+c_0|}{v}) = -\frac{\pi}{2} +O\big(\frac{|\sigma|}{v}\big),
\]
hence
\begin{equation}
\label{L3-1/3}
\frac{1}{\Gamma(1-w)} \ll |1-w|^{u-\frac12} e^{-v(-\frac{\pi}{2} + O(\frac{|\sigma|}{v}))} A^{|\sigma|}
\ll e^{\frac{\pi}{2}v} (v+|\sigma|)^{|\sigma|+c} A^{|\sigma|}.
\end{equation}
The result follows from \eqref{L3-1/1}-\eqref{L3-1/3} since $d_F=2$.  \fine

\bigskip
{\bf Lemma 3.2.}  {\sl For $w\in\LL_\infty(s)$ and $s\in\CC$ we have
\[
S(s,w) = -ie(-\xi_F/4)e^{-\pi i s}\big(1+O(A^{|s|}e^{-cv})\big) \hskip2cm v\to+\infty
\]
with suitable constants} $A,c>0$.

\bigskip
{\sl Proof.} Since $d_F=2$ we have
\[
\begin{split}
S(s,w) &= 2^{r-1}(2i)^{1-r} \frac{\prod_{j=1}^r \big(e^{i\pi(\lambda_j(s+w)+\mu_j)} - e^{-i\pi(\lambda_j(s+w)+\mu_j)}\big)}{e^{i\pi w} - e^{-i\pi w}} \\
&= (-i)^{1-r} e^{i\pi w} \prod_{j=1}^r e^{-i\pi(\lambda_j(s+w)+\mu_j)}  \big(1+O(A^{|s|}e^{-cv})\big) \\
&=-ie^{-i\pi s -i\pi \sum_{j=1}^r(\mu_j-\frac12)}  \big(1+O(A^{|s|}e^{-cv})\big)
\end{split}
\]
and the result follows. \fine

\bigskip
{\bf Lemma 3.3.} {\sl For $s\in\CC$ we have}
\[
\begin{split}
R_\nu(s) &= B_{\nu+1}(1-2s-i\theta_F) + B_{\nu+1}(1) \\
&- \sum_{j=1}^r \frac{B_{\nu+1}(\lambda_j(1-s)+\bar{\mu}_j) + B_{\nu+1}(1-\lambda_js-\mu_j)}{\lambda_j^\nu}.
\end{split}
\]

\bigskip
{\sl Proof.} We use the following properties of the Bernoulli polynomials (see Section 1.13 of Bateman's Project \cite{EMOT/1953}):
\begin{equation}
\label{L3-3/1}
B_n(1-x) = (-1)^nB_n(x)
\end{equation}
\begin{equation}
\label{L3-3/2}
B_n(x+y) = \sum_{k=0}^n {n\choose k} B_k(x)y^{n-k}.
\end{equation}
By \eqref{L3-3/1} and \eqref{L3-3/2} we have
\[
\begin{split}
B_{\nu+1}(1-\lambda_js-\mu_j) &= (-1)^{\nu+1} B_{\nu+1}(\lambda_js+\mu_j) \\
&= (-1)^{\nu+1} \sum_{k=0}^{\nu+1} {\nu+1 \choose k} B_k(\mu_j)\lambda_j^{\nu+1-k} s^{\nu+1-k}.
\end{split}
\]
Similarly
\[
B_{\nu+1}(\bar{\mu}_j +\lambda_j(1-s)) = \sum_{k=0}^{\nu+1} {\nu+1\choose k} B_k(\bar{\mu}_j)\lambda_j^{\nu+1-k}(1-s)^{\nu+1-k}.
\]
Therefore
\[
\begin{split}
\sum_{j=1}^r &\frac{B_{\nu+1}(\lambda_j(1-s)+\bar{\mu}_j) + B_{\nu+1}(1-\lambda_js-\mu_j)}{\lambda_j^\nu} \\
&= \sum_{j=1}^r \sum_{k=0}^{\nu+1} {\nu+1\choose k} \big( \frac{B_k(\bar{\mu}_j)}{\lambda_j^{k-1}} (1-s)^{\nu+1-k} + (-1)^{\nu+1}\frac{B_k(\mu_j)}{\lambda_j^{k-1}} s^{\nu+1-k}\big) \\
&= -\frac12 \sum_{k=0}^{\nu+1} {\nu+1\choose k} \big((-1)^\nu H_F(k) s^{\nu+1-k} - \overline{H_F(k)}(1-s)^{\nu+1-k}\big)
\end{split}
\]
and the result follows from the definition of the $R_\nu(s)$'s. \fine

\bigskip
{\bf Lemma 3.4.} {\sl For $w\in\LL_\infty(s)$ and $s\in\CC$ we have
\[
\Gamma(1-2s-w-i\theta_F) \ll e^{-\frac{\pi}{2}v} A^{|s|} v^{-\sigma+2c_0} \hskip2cm v\to+\infty
\]
with a suitable constant} $A>0$.

\bigskip
{\sl Proof.} Since
\[
\arg(1-2s-w-i\theta_F) = -\frac{\pi}{2} + O(\frac{|s|}{v})
\]
we have 
\[
\begin{split}
\log&|\Gamma(1-2s-w-i\theta_F)| \\
&= (\frac12-\sigma+c_0)\log|1-2s-w-i\theta_F| + (v+2t+\theta_F)\arg(1-2s-w-i\theta_F) \\
&= (\frac12-\sigma+c_0)\big(\log v +  O(\frac{|s|}{v})\big) -\frac{\pi}{2} v + O(|s|),
\end{split}
\]
and the result follows since $|s|\ll v$. \fine

\bigskip
{\bf Lemma 3.5.} {\sl For $\nu\geq 1$, $0<D<2\pi$ and $s\in\CC$ we have $B_\nu(s) \ll \nu! \frac{e^{D|s|}}{(2\pi-D)D^\nu}$; in particular}
\[
B_\nu(s) \ll e^{|s|} \nu!.
\]

\bigskip
{\sl Proof.} By definition we have for $|z|<2\pi$
\[
\frac{ze^{zs}}{e^z-1} = \sum_{\nu=0}^\infty \frac{B_\nu(s)}{\nu!} z^\nu.
\]
Hence, since the function $\frac{z}{e^z-1}$ is holomorphic on $|z|\leq 3\pi$ (say) apart from two simple poles at $z=\pm2\pi i$, we have
\[
\big|\frac{B_\nu(s)}{\nu!}\big| = \big|\frac{1}{2\pi i} \int_{|z|=D}  \frac{ze^{zs}}{e^z-1} \frac{\d z}{z^{\nu+1}}\big| \ll \frac{e^{D|s|}}{D^{\nu+1}} \int_{|z|=D} \big|\frac{z}{e^z-1}\big|\d z \ll \frac{e^{D|s|}}{(2\pi-D)D^\nu}
\]
and the result follows. \fine

\bigskip
{\bf Lemma 3.6.}  {\sl For $1\leq N \leq \frac{3}{2}x$ and $x>0$ we have}
\[
\Phi_N(x):= \sum_{m=1}^N\frac{x^m}{m!} \leq \frac{(2x)^N}{N!}.
\]

\bigskip
{\sl Proof.} This is Lemma 7 of \cite{Ka-Pe/USF}, after deleting the unneeded $-1$ in its statement. \fine

\bigskip
{\bf Lemma 3.7.}  {\sl For $F\in\S^\sharp$ and $\nu\geq 1$ we have
\[
H_F(\nu) \ll c^\nu \nu!
\]
with a suitable} $c>0$.

\bigskip
{\sl Proof.} By definition we have
\[
H_F(\nu) = 2\sum_{j=1}^r \frac{B_\nu(\mu_j)}{\lambda_j^{\nu-1}} \ll c^\nu \max_j|B_\nu(\mu_j)|,
\]
and the result follows thanks to Lemma 3.5. \fine

\bigskip
{\bf Lemma 3.8.}  {\sl Let $c\geq 1$. For $s\in\CC$ and $1\leq \nu\leq c(|s|+1)$ we have
\[
R_\nu(s) \ll \big(c'(|s|+1)\big)^{\nu+1}
\]
with a suitable $c'>0$ (depending also on $c$).}

\bigskip
{\sl Proof.} Suppose first that $|s|\leq 1$; then $\nu$ is bounded and hence $R_\nu(s)$ is also bounded, and the result follows in this case. Let now $|s|>1$. Since $B_{\nu+1}(1) = B_{\nu+1}$ (see Section 1.13 of  \cite{EMOT/1953}) and ${\nu+1\choose k}\leq 2^{\nu+1}$ we have
\begin{equation}
\label{L3-8/1}
R_\nu(s) \ll |B_{\nu+1}(1-2s-i\theta_F)| +|B_{\nu+1}| + \big(c(|s|+1)\big)^{\nu+1} \sum_{k=0}^{\nu+1} |H_F(k)|s|^{-k}.
\end{equation}
Since $B_k=B_k(0)$, from Lemma 3.5 we have for every $k\in\NN$
\[
B_k\ll k!.
\]
Moreover, from \eqref{L3-3/2} we have
\[
B_\nu(x) = \sum_{k=0}^\nu{\nu\choose k} B_kx^{\nu-k}
\]
and hence by Lemma 3.6 for $\nu\leq c(|s|+1)$ we get
\begin{equation}
\label{L3-8/2}
\begin{split}
B_{\nu+1}(1-2s+i\theta_F) &\ll \sum_{k=0}^{\nu+1} {\nu+1\choose k} k! \big(c'(|s|+1)\big)^{\nu+1-k}  \\
&\ll (\nu+1)!\big(\Phi_{\nu+1}(c'(|s|+1))+1\big) \ll \big(c'(|s|+1)\big)^{\nu+1}.
\end{split}
\end{equation}
Further
\begin{equation}
\label{L3-8/3}
B_{\nu+1} \ll (\nu+1)! \ll  \big(c'(|s|+1)\big)^{\nu+1},
\end{equation}
and finally thanks to Lemma 3.7 we obtain
\begin{equation}
\label{L3-8/4}
 \sum_{k=0}^{\nu+1} |H_F(k)|s|^{-k} \ll \tilde{c}^{\nu+1} \sum_{k=0}^{\nu+1} \frac{k!}{|s|^k} \ll \tilde{c}^{\nu+1} \sum_{k=0}^{\nu+1} \big(\frac{k}{|s|}\big)^k \ll \tilde{c}^{\nu+1} \sum_{k=0}^{\nu+1} \big(\frac{\nu+1}{|s|}\big)^k \ll \tilde{c}^{\nu+1},
\end{equation}
and the result follows from \eqref{L3-8/1}-\eqref{L3-8/4}. \fine

\bigskip
{\bf Lemma 3.9.}  {\sl Let $c\geq 1$. For $w\in\LL_\infty(s)$, $s\in\CC$ and $1\leq N\leq c(|s|+1)$ we have
\[
\sum_{\nu=1}^N \frac{R_\nu(s)}{\nu(\nu+1)} \frac{1}{w^\nu} \leq c'(|s|+1).
\]
with a suitable $c'>0$ (depending also on $c$).}

\bigskip
{\sl Proof.} By Lemma 3.8 and choosing $c_1$ and $c_2$ in the definition of $\LL_\infty(s)$ sufficiently large we have
\[
\sum_{\nu=1}^N \frac{R_\nu(s)}{\nu(\nu+1)} \frac{1}{w^\nu} \ll \sum_{\nu=1}^N \frac{\big(c'(|s|+1)\big)^{\nu+1}}{\nu(\nu+1)} \frac{1}{|w|^\nu} \ll c'(|s|+1)\sum _{\nu=1}^N \frac{1}{\nu^2},
\]
and the result follows. \fine

\bigskip
{\bf Lemma 3.10.}  {\sl For $w\in\CC$ and integers $m,M$ with $1\leq m+1\leq M<|w|$ we have
\[
\frac{1}{w} = \frac{(-1)^{m+1}}{m!} \sum_{\ell=m+1}^M \frac{(-1)^\ell (\ell-1)!}{(w-(m+1))\cdots(w-\ell)} + r_{m,M}(w)
\]
with}
\[
|r_{m,M}(w)| \leq \frac{M!}{m!|w(w-(m+1))\cdots(w-M)|}.
\]

\bigskip
{\sl Proof.} We proceed by induction on $m$. For $m=0$ one easily verifies the identity
\[
\frac{1}{w} = -\sum_{\ell=1}^M \frac{(-1)^\ell (\ell-1)!}{(w-1)\cdots(w-\ell)} + \frac{(-1)^M M!}{w(w-1)\cdots(w-M)}
\]
by induction on $M$, and the lemma follows in this case. Assume now that the result holds true for a certain $m\geq 0$ and let $M\geq m+2$. We have
\[
\frac{-(m+1)}{w(w-(m+1))} = \frac{1}{w} - \frac{1}{w-(m+1)} = \frac{(-1)^{m+1}}{m!} \sum_{\ell=m+2}^M \frac{(-1)^\ell (\ell-1)!}{(w-(m+1))\cdots(w-\ell)} + r_{m,M}(w),
\]
hence multiplying by $-\frac{w-(m+1)}{m+1}$ we obtain
\[
\frac{1}{w} = \frac{(-1)^{m+2}}{(m+1)!} \sum_{\ell=m+2}^M \frac{(-1)^\ell (\ell-1)!}{(w-(m+2))\cdots(w-\ell)} + r_{m+1,M}(w)
\]
with
\[
 r_{m+1,M}(w) = -\frac{w-(m+1)}{m+1}r_{m,M}(w).
\]
Therefore
\[
\begin{split}
|r_{m+1,M}(w)| &\leq \frac{|w-(m+1)|}{m+1} \frac{M!}{m!|w(w-(m+1))\cdots(w-M)|} \\
&= \frac{M!}{(m+1)!|w(w-(m+2))\cdots(w-M)|},
\end{split}
\]
and the result follows by induction. \fine

\bigskip
Now we obtain a recursive relation with respect to $\mu$ for the coefficients $C_{\mu,\ell}$ in \eqref{3-2}. Since we did not yet establish \eqref{3-2}, the relations in Lemma 3.11 below are in a first instance to be considered as relations between coefficients of formal expansions. However, \eqref{3-2} is then established (in a precise form) in Lemma 3.13 below.

\bigskip
{\bf Lemma 3.11.}  {\sl Formally, for $\ell\geq \mu+1\geq 2$ we have}
\[
C_{\mu+1,\ell} = (-1)^{\ell-1} (\ell-1)! \sum_{m=\mu}^{\ell-1} \frac{(-1)^m C_{\mu,m}}{m!}.
\]

\bigskip
{\sl Proof.} As remarked above, we proceed using formal expansions. Using Lemma 3.10 with each $m\geq \mu$ we have
\[
\begin{split}
\frac{1}{w^{\mu+1}} = \frac{1}{w^\mu}\frac{1}{w} &= \sum_{m=\mu}^\infty \frac{C_{\mu,m}}{(w-1)\cdots(w-m)} \big(\frac{(-1)^{m+1}}{m!} \sum_{\ell=m+1}^\infty \frac{(-1)^\ell (\ell-1)!}{(w-(m+1))\cdots(w-\ell)}\big) \\
& = \sum_{\ell = \mu+1}^\infty \big((-1)^{\ell-1} (\ell-1)! \sum_{m=\mu}^{\ell-1} \frac{(-1)^m C_{\mu,m}}{m!} \big) \frac{1}{(w-1)\cdots(w-\ell)}
\end{split}
\]
and the result follows. \fine

\bigskip
{\bf Lemma 3.12.}  {\sl For $\ell \geq 1$ we have 
\[
C_{1,\ell} = (-1)^{\ell-1}(\ell-1)!.
\]
Moreover,  if the $C_{\mu,\ell}$'s satisfy the recurrence in Lemma $3.11$, then for $\ell \geq \mu\geq 2$
\[
C_{\mu,\ell} = (-1)^{\ell-\mu}(\ell-1)! \sum_{1\leq k_1<k_2<\cdots<k_{\mu-1}\leq \ell-1} \frac{1}{k_1k_2\cdots k_{\mu-1}}
\]
and for} $\ell\geq \mu \geq 1$
\[
|C_{\mu,\ell}| \leq \frac{(\ell-1)!}{(\mu-1)!}{\ell-1\choose \mu-1}.
\]

\bigskip
{\sl Proof.} The first assertion, for $\ell\geq \mu=1$, follows from Lemma 3.10 with $m=0$. To prove the second assertion we argue by induction on $\mu\geq 2$. For $\ell\geq \mu=2$ we have from Lemma 3.11 that
\[
C_{2,\ell} = (-1)^{\ell-1} (\ell-1)! \sum_{m=1}^{\ell-1} \frac{(-1)^m C_{1,m}}{m!} = (-1)^{\ell-2} (\ell-1)! \sum_{m=1}^{\ell-1} \frac{1}{m},
\]
and the result follows in this case. For $\ell\geq \mu+1$, inserting the inductive hypothesis in the expression given by Lemma 3.11 we obtain
\[
\begin{split}
C_{\mu+1,\ell} &= (-1)^{\ell-1} (\ell-1)! \sum_{m=\mu}^{\ell-1} \frac{(-1)^m}{m!}
(-1)^{m-\mu}(m-1)! \sum_{1\leq k_1<k_2<\cdots<k_{\mu-1}\leq m-1} \frac{1}{k_1k_2\cdots k_{\mu-1}} \\
&= (-1)^{\ell-(\mu+1)} (\ell-1)! \sum_{m=\mu}^{\ell-1} \ \sum_{1\leq k_1<k_2<\cdots<k_{\mu-1}\leq m-1} \frac{1}{mk_1k_2\cdots k_{\mu-1}} \\
&= (-1)^{\ell-(\mu+1)} (\ell-1)! \sum_{1\leq k_1<k_2<\cdots<k_\mu\leq \ell-1} \frac{1}{k_1k_2\cdots k_\mu}
\end{split}
\]
and the second assertion follows by induction. The third assertion follows easily from the previous ones since $k_1k_2\cdots k_{\mu-1} \geq (\mu-1)!$. \fine

\bigskip
{\bf Lemma 3.13.}  {\sl For $w\in\CC$ and integers $\mu,M$ with $1\leq \mu\leq M\leq |w|/2$ we have
\[
\frac{1}{w^\mu} = \sum_{\ell =\mu}^M \frac{C_{\mu,\ell}}{(w-1)\cdots(w-\ell)} + R_{\mu,M}(w)
\]
with}
\[
R_{\mu,M}(w) \ll \frac{2^MM!}{(\mu-1)!} \frac{1}{|w(w-1)\cdots(w-M)|}.
\]

\bigskip
{\sl Proof.} By induction on $\mu$. For $\mu=1$ the assertion follows at once from Lemma 3.10 with $m=0$ and Lemma 3.12. From Lemma 3.10, arguing similarly as in Lemma 3.11, we obtain
\[
\begin{split}
\frac{1}{w^{\mu+1}} = \frac{1}{w^\mu} \frac{1}{w} &= \sum_{m=\mu}^M \frac{C_{\mu,m}}{(w-1)\cdots(w-m)} \big(\frac{(-1)^{m+1}}{m!} \sum_{\ell=m+1}^M \frac{(-1)^\ell (\ell-1)!}{(w-(m+1))\cdots(w-\ell)}\big) \\
&+ R_{\mu+1,M}(w)
\end{split}
\]
with
\[
R_{\mu+1,M}(w) =  \sum_{m=\mu}^M \frac{C_{\mu,m}}{(w-1)\cdots(w-m)} r_{m,M}(w) + \frac{1}{w} R_{\mu,M}(w).
\]
Similarly as in Lemma 3.11, from the inductive hypothesis we obtain
\[
\frac{1}{w^{\mu+1}} = \sum_{\ell=\mu+1}^M \frac{C_{\mu+1,\ell}}{(w-1)\cdots(w-\ell)} + R_{\mu+1,M}(w).
\]
Moreover, from the inductive hypothesis coupled with Lemmas 3.10 and 3.12 we get
\[
\begin{split}
R_{\mu+1,M}(w) &\ll \sum_{m=\mu}^M \frac{(m-1)!}{(\mu-1)!} {m-1\choose \mu-1} \frac{M!}{m!|w(w-1)\cdots(w-M)|} \\
&+ \frac{1}{|w|} \frac{2^MM!}{(\mu-1)!} \frac{1}{|w(w-1)\cdots(w-M)|} \\
&\leq \frac{M!}{(\mu-1)!|w(w-1)\cdots(w-M)|} \big(\sum_{m=\mu}^M \frac{{m-1\choose \mu-1}}{m} + \frac{2^M}{|w|}\big) \\
&\leq \frac{2^MM!}{\mu!|w(w-1)\cdots(w-M)|} \big(2^{-M}\sum_{m=\mu}^M {m-1\choose \mu-1} + \frac{\mu}{|w|}\big) \\
&\leq \frac{2^MM!}{\mu!|w(w-1)\cdots(w-M)|},
\end{split}
\]
since $\mu/|w| \leq 1/2$ and
\[
2^{-M}\sum_{m=\mu}^M {m-1\choose \mu-1} = 2^{-M}\sum_{m=\mu}^M {m-1\choose m-\mu} \leq 2^{-M}\sum_{m=\mu}^M {M-1\choose m-\mu} \leq 2^{-M}2^{M-1} = \frac12.
\]
This ends the proof of the lemma. \fine

\bigskip
Now we turn to the asymptotic expansion in \eqref{3-3}. The next lemma holds for any $\theta_F\in\CC$.

\bigskip
{\bf Lemma 3.14.}  {\sl Let $s\in\CC$ and $1\leq \mu\leq\nu\leq N$ be integers. Then
\[
A_{\mu,\nu}(s) = \sum_{k=0}^{\nu-\mu} {-\mu\choose k} C_{\mu+k,\nu}(2s-1+i\theta_F)^k.
\]
Moreover, if $1 \leq c\leq c_2/3$, $w\in\LL_{-\infty}^*(s)$ and $1\leq N \leq |\sigma| +c$ we have
\[
\frac{1}{(w+2s-1+i\theta_F)^\mu} = \sum_{\nu=\mu}^N \frac{A_{\mu,\nu}(s)}{(w-1)\cdots(w-\nu)} + O\big(\frac{A^{|s|}(|s|^{|\sigma|}+1)}{(\mu-1)!} \frac{1}{|w(w-1)\cdots(w-N)|}\big)
\]
with a suitable $A>0$. The constant in the $O$-symbol may depend also on $c$.}

\bigskip
{\sl Proof.} We first observe that for $|\theta|<1/2$ and integer $P\geq 0$
\begin{equation}
\label{L3-14/1}
\frac{1}{(1+\theta)^\mu} = \sum_{k=0}^P {-\mu\choose k} \theta^k + O\big(\frac{2^{P+\mu}|\theta|^{P+1}}{1-2|\theta|}\big).
\end{equation}
Indeed we have, since $\big|{-\mu\choose k}\big| = {\mu+k-1\choose k}$, that
\[
\sum_{k=P+1}^\infty \big|{-\mu\choose k}\big| |\theta|^k = \sum_{k=P+1}^\infty {\mu+k-1\choose k} |\theta|^k \leq 2^{\mu-1} \sum_{k=P+1}^\infty |2\theta|^k = \frac{2^{P+\mu}|\theta|^{P+1}}{1-2|\theta|},
\]
and \eqref{L3-14/1} follows. For brevity we write $\eta = 2s-1+i\theta_F$ and $\theta= \eta/w$. Since $1-2|\theta| \geq 1/2$ for $w\in\LL_{-\infty}^*(s)$ provided $c_1$ and $c_2$ are sufficiently large, given an integer $P\geq 0$ from \eqref{L3-14/1} we obtain
\begin{equation}
\label{L3-14/2}
\begin{split}
\frac{1}{(w+\eta)^\mu} &= \frac{1}{w^\mu(1+\theta)^\mu} = \frac{1}{w^\mu}\left(\sum_{k=0}^P {-\mu\choose k} \theta^k + O(\frac{A^{P+\mu}(|s|+1)^{P+1}}{|w|^{P+1}})\right) \\
&= \sum_{k=0}^P {-\mu\choose k} \frac{\eta^k}{w^{\mu+k}} + O\big( \frac{A^{P+\mu}(|s|+1)^{P+\mu+1}}{(|s|+1)^\mu|w|^{P+\mu+1}}\big),
\end{split}
\end{equation}
where the last $O$-term is written in a suitable way for the estimate below. Choose $P=N-\mu$. Since $|w|\geq 2N$ for $w\in\LL_{-\infty}^*$ we have
\[
|w|^{N+1} \geq 2^{-N}|w(w-1)\cdots(w-N)|.
\]
Moreover, since $\mu\leq N\leq |\sigma| +c$ we have $(|s|+1)^\mu \gg (\mu-1)!$, and the error term in \eqref{L3-14/2} is
\begin{equation}
\label{L3-14/3}
\ll \frac{A^{|s|}(|s|^{|\sigma|}+1)}{(\mu-1)!} \frac{1}{|w(w-1)\cdots(w-N)|}.
\end{equation}
Now we apply Lemma 3.13 with $M=N$, thus getting that the main term in \eqref{L3-14/2} is
\begin{equation}
\label{L3-14/4}
\begin{split}
&\sum_{k=0}^{N-\mu}{-\mu\choose k} \eta^k \sum_{\nu=\mu+k}^N \frac{C_{\mu+k,\nu}}{(w-1)\cdots(w-\nu)} \\
&+ O\big(\sum_{k=0}^{N-\mu}\big|{-\mu\choose k}\big| |\eta|^k \frac{2^N N!}{(\mu+k-1)!} \frac{1}{|w(w-1)\cdots(w-N)|}\big) \\
&= \sum_{\nu=\mu}^N \frac{\sum_{k=0}^{\nu-\mu} {-\mu\choose k} C_{\mu+k,\nu}\eta^k}{(w-1)\cdots(w-\nu)} +O\big(\frac{A^{|s|}(|s|^{|\sigma|}+1)}{(\mu-1)!} \sum_{k=0}^N \frac{|\eta|^k}{k!} \frac{1}{|w(w-1)\cdots(w-N)|}\big),
\end{split} 
\end{equation}
since
\[
2^NN! \leq 2^N N^N \ll A^{|s|}(|s|^{|\sigma|}+1).
\]
Moreover
\[
\sum_{k=0}^N \frac{|\eta|^k}{k!} \ll e^{|\eta|} \ll A^{|s|},
\]
therefore the error terms in \eqref{L3-14/3} and \eqref{L3-14/4} are of the required size. Hence both the expression for the $A_{\mu,\nu}(s)$ and the asymptotic expansion of $1/(w+\eta)^\mu$ follow from \eqref{L3-14/2} and \eqref{L3-14/4}, and the lemma is proved. \fine

\bigskip
{\bf Lemma 3.15.}  {\sl Let $s\in\CC$ and $1 \leq c \leq \min(c_1,c_2)/10$. For $w\in\LL_{-\infty}^*(s)$ and integer $1\leq N\leq |\sigma|+c$ we have
\[
\begin{split}
\exp\big(\sum_{\nu=1}^N  \frac{(-1)^\nu R_\nu(s)}{\nu(\nu+1)} \frac{1}{(w+2s-1+i\theta_F)^\nu}\big) &= \sum_{\nu=0}^N \frac{Q_\nu(s)}{(w-1)\cdots(w-\nu)} \\ 
&+ O\big(\frac{A^{|s|}(|s|^{|\sigma|}+1)}{|w(w-1)\cdots(w-N)|}\big)
\end{split}
\]
where the $Q_\nu(s)$'s are polynomials with $Q_0(s)=1$ identically, $A>0$ is suitable and the constant in the $O$-symbol may depend also on $c$. Moreover, for every $\nu\geq 1$ we have}
\begin{equation}
\label{L3-15/0}
Q_\nu(s) = \sum_{\mu=1}^\nu V_\mu(s)A_{\mu,\nu}(s).
\end{equation}

\bigskip
{\sl Proof.} Writing again $\eta = 2s-1+i\theta_F$, from the power series expansion of the exponential function we have
\begin{equation}
\label{L3-15/1}
\exp\big(\sum_{\nu=1}^N  \frac{(-1)^\nu R_\nu(s)}{\nu(\nu+1)} \frac{1}{(w+\eta)^\nu}\big) = 1 + \sum_{\mu=1}^\infty \frac{V_{\mu,N}(s)}{(w+\eta)^\mu},
\end{equation}
where
\[
V_{\mu,N}(s) = (-1)^\mu \sum_{m=1}^\mu \frac{1}{m!}  \sum_{\substack{1\leq \nu_1\leq N,...,1\leq \nu_m\leq N\\ \nu_1+...+\nu_m=\mu}} \prod_{j=1}^m\frac{R_{\nu_j}(s)}{\nu_j(\nu_j+1)}.
\]
Recalling definition \eqref{3-1} we see that $V_{\mu,N}(s) = V_{\mu}(s)$ for $\mu\leq N$, while for every $\mu\geq 1$ thanks to Lemma 3.8 we obtain
\begin{equation}
\label{L3-15/2}
\begin{split}
|V_{\mu,N}(s)| &\leq \sum_{m=1}^\mu \frac{1}{m!}  \sum_{ \nu_1+...+\nu_m=\mu} \prod_{j=1}^m\frac{(c(|s|+1))^{\nu_j+1}}{\nu_j(\nu_j+1)} \\
&\ll (c(|s|+1))^\mu \sum_{m=1}^\mu \frac{(c(|s|+1))^m}{m!} \ll A^{|s|} (c(|s|+1))^\mu.
\end{split}
\end{equation}
Hence, from \eqref{L3-15/2}, for $w\in\LL_{-\infty}^*$ and $N\leq |\sigma|+c$ we have
\begin{equation}
\label{L3-15/3}
\begin{split}
\sum_{\mu=N+1}^\infty \frac{|V_{\mu,N}(s)|}{|w+\eta|^\mu} &\ll A^{|s|} \sum_{\mu=N+1}^\infty \frac{(c'(|s|+1))^\mu}{|w+\eta|^\mu} \ll A^{|s|} \big(\frac{c'(|s|+1)}{|w+\eta|}\big)^{N+1} \\
&\ll A^{|s|} \frac{|s|^{|\sigma|}+1}{|w|^{N+1}} \ll \frac{A^{|s|}(|s|^{|\sigma|}+1)}{|w(w-1)\cdots(w-N)|}.
\end{split}
\end{equation}
Finally, using Lemma 3.14 and again  \eqref{L3-15/2} we get
\begin{equation}
\label{L3-15/4}
\begin{split}
\sum_{\mu=1}^N \frac{V_{\mu}(s)}{(w+\eta)^\mu} &= \sum_{\mu=1}^N V_\mu(s) \sum_{\nu=\mu}^N \frac{A_{\mu,\nu}(s)}{(w-1)\cdots(w-\nu)} \\
&+ O\big(\frac{A^{|s|}(|s|^{|\sigma|}+1)}{|w(w-1)\cdots(w-N)|} \sum_{\mu=1}^N \frac{(c'(|s|+1))^\mu}{(\mu-1)!}\big) \\
&= \sum_{\nu=1}^N \frac{Q_\nu(s)}{(w-1)\cdots(w-\nu)} + O\big(\frac{A^{|s|}(|s|^{|\sigma|}+1)}{|w(w-1)\cdots(w-N)|}\big),
\end{split}
\end{equation}
and the result follows from \eqref{L3-15/1}, \eqref{L3-15/3}, \eqref{L3-15/4} and the fact that the $V_\mu(s)$'s and the $A_{\mu,\nu}(s)$'s are polynomials. \fine

\bigskip
{\bf Lemma 3.16.}  {\sl Let $s\in\CC$ and $c\geq 1$. For $1\leq \mu \leq c(|s|+1)$ we have
\[
|V_\mu(s)| \leq \frac{(c'(|s|+1))^{2\mu}}{\mu!}
\]
with suitable $c'=c'(c)>0$}.

\bigskip
{\sl Proof.} As in the proof of Lemma 3.15 (see \eqref{L3-15/2} and recall that $V_\mu(s) = V_{\mu,\mu}(s)$) we have, with the notation of Lemma 3.6, that
\[
|V_\mu(s)| \ll \big(c^*(|s|+1)\big)^\mu \sum_{m=1}^\mu \frac{(c^*(|s|+1))^m}{m!} = \big(c^*(|s|+1)\big)^\mu \Phi_\mu(c^*(|s|+1)).
\]
Here $c^*$ is chosen so large that Lemma 3.6 is applicable, hence the result follows from Lemma 3.6. \fine

\bigskip
{\bf Lemma 3.17.}  {\sl Let $x>0$, $M$ be a positive integer and
\[
\Psi_M(x) = \sum_{\mu=1}^M \frac{x^{2\mu}}{(\mu!)^2}.
\]
Then for $1\leq M \leq \sqrt{2}x$ we have}
\[
\Psi_M(x) \leq \frac{(2x)^{2M}}{(M!)^2}.
\]

\bigskip
{\sl Proof.} By induction on $M$. It is trivial for $M=1$, and assume it holds for an integer $M$. Then since $1\leq M+1 \leq \sqrt{2}x$ we have
\[
\begin{split}
\Psi_{M+1}(x) &= \Psi_M(x) + \frac{x^{2(M+1)}}{((M+1)!)^2} \leq \frac{(2x)^{2M}}{(M!)^2} + \frac{x^{2(M+1)}}{((M+1)!)^2} \\
&= \frac{(2x)^{2(M+1)}}{((M+1)!)^2} \big\{\frac{(M+1)^2}{(2x)^2} + \frac{1}{4^{M+1}}\big\} \leq \frac{9}{16} \frac{(2x)^{2(M+1)}}{((M+1)!)^2},
\end{split} 
\]
and the result follows. \fine

\bigskip
{\bf Lemma 3.18.}  {\sl Let $s\in\CC$ and $c\geq 1$. For $1\leq \nu\leq |s|+c$ we have
\[
Q_\nu(s) \ll \frac{(c'(|s|+1))^{2\nu}}{\nu!}
\]
with suitable} $c'=c'(c)>0$.

\bigskip
{\sl Proof.} From Lemmas 3.15, 3.14, 3.16 and 3.12 and the fact that $\big|{-\mu\choose k}\big| = {\mu+k-1\choose k}$ we have
\[
\begin{split}
Q_\nu(s)& \ll \sum_{\mu=1}^\nu |V_\mu(s)| |A_{\mu,\nu}(s)| \ll \sum_{\mu=1}^\nu \sum_{k=0}^{\nu-\mu} {\mu+k-1 \choose k} |C_{\mu+k,\nu}| (c^*(|s|+1))^k \frac{(c^*(|s|+1))^{2\mu}}{\mu!} \\
&\ll \sum_{\mu=1}^\nu \sum_{k=0}^{\nu-\mu} {\mu+k-1 \choose k}  {\nu-1\choose \mu+k-1} \frac{(\nu-1)!}{(\mu+k-1)!} (c^*(|s|+1))^k \frac{(c^*(|s|+1))^{2\mu}}{\mu!} \\
&\ll 4^\nu \nu! \sum_{k=0}^{\nu-1} (c^*(|s|+1))^k \sum_{\mu=1}^{\nu-k} \frac{(c^*(|s|+1))^{2\mu}}{\mu!(\mu+k)!}
\end{split}
\]
since
\[
{\mu+k-1\choose k} \leq 2^{\mu+k-1} \leq 2^\nu, \qquad {\nu-1\choose \mu+k-1} \leq 2^\nu, \qquad \frac{(\nu-1)!}{(\mu+k-1)!}  \leq \frac{\nu!}{(\mu+k)!}.
\]
Moreover, using $(\mu+k)! \geq \mu!k!$ we see that the last expression is (see Lemma 3.17)
\[
\ll  4^\nu \nu! \sum_{k=0}^{\nu-1} \frac{(c^*(|s|+1))^k}{k!} \Psi_{\nu-k}(c^*(|s|+1)).
\]
Hence by Lemmas 3.17 and 3.6 we get, provided $c^*$ is large enough with respect to $c$,
\[
\begin{split}
Q_\nu(s) &\ll  4^\nu \nu! \sum_{k=0}^{\nu-1} \frac{(c^*(|s|+1))^k}{k!} \frac{ (c^*(|s|+1))^{2(\nu-k)}}{((\nu-k)!)^2} \ll (c^*(|s|+1))^\nu  \sum_{k=0}^{\nu-1} {\nu\choose k} \frac{(c^*(|s|+1))^{\nu-k}}{(\nu-k)!} \\
&= (c^*(|s|+1))^\nu  \Phi_\nu(c^*(|s|+1)) \ll \frac{(c'(|s|+1))^{2\nu}}{\nu!},
\end{split}
\]
and the result follows. \fine

\bigskip
{\bf Lemma 3.19.}  {\sl Let $K>0$ be an integer and $c\geq 1$. Then for $|s| \leq 2K$ and $1\leq\nu\leq K+c$ we have
\begin{equation}
\label{L3-19/1}
R_\nu(s) \ll A^K K^{\nu},
\end{equation}
\begin{equation}
\label{L3-19/2}
V_\nu(s) \ll A^K K^\nu
\end{equation}
\begin{equation}
\label{L3-19/3}
Q_\nu(s) \ll (AK)^K
\end{equation}
with a suitable constant $A=A(c)>0$. Moreover,} $\deg Q_\nu(s) = 2\nu$.

\bigskip
{\sl Proof.} We start with \eqref{L3-19/1}. If $\nu\leq |s|$ and $|s| \leq 2K$ then Lemma 3.8 gives
\[
R_\nu(s) \ll (c(|s|+1))^{\nu+1} \ll A^K K^{\nu}.
\]
Suppose now that $|s| \leq \nu \leq K+c$. Then by Lemmas 3.5 and 3.7 we get
\[
\begin{split}
R_\nu(s) &\ll |B_{\nu+1}(1-2s-i\theta_F)| + |B_{\nu+1}| + 2^{\nu} \sum_{k=0}^{\nu+1} |H_F(k)| (|s|+1)^{\nu+1-k} \\
&\ll e^{2|s|} (\nu+1)! + A^{\nu} \sum_{k=0}^{\nu+1} k! (|s|+1)^{\nu+1-k} \\
&\ll e^{2|s|} (\nu+1)! + A^{\nu}(|s|+1)^{\nu+1}\sum_{k=0}^{\nu+1}\big(\frac{k}{|s|+1}\big)^k \\
&\ll A^K K^{\nu} +A^\nu(\nu+1)^\nu \ll A^K K^{\nu}
\end{split}
\]
and \eqref{L3-19/1} follows. To prove  \eqref{L3-19/2} we observe that from the definition of the $V_\nu(s)$'s and \eqref{L3-19/1} we have
\[
V_\nu(s) \ll \sum_{m=1}^\nu \frac{1}{m!} \sum_{\nu_1+\cdots+\nu_m=\nu} \prod_{j=1}^m\frac{A^K K^{\nu_j+1}}{\nu_j(\nu_j+1)} \ll A^K K^\nu\sum_{m=1}^\nu \frac{K^m}{m!}\big(\sum_{\ell=1}^\infty \frac{1}{\ell(\ell+1)}\big)^m \ll  A^K K^\nu,
\]
thus proving \eqref{L3-19/2}. Finally, starting as in Lemma 3.18 and using \eqref{L3-19/2} instead of Lemma 3.16  we get
\[
\begin{split}
Q_\nu(s) &\ll A^K \sum_{\mu=1}^\nu \sum_{k=0}^{\nu-\mu} {\mu+k-1 \choose k} |C_{\mu+k,\nu}| (AK)^{\mu+k} \\
&\ll A^K 4^\nu \nu! \sum_{\mu=1}^\nu \sum_{k=0}^{\nu-\mu}  \frac{(AK)^{\mu+k}}{(\mu+k)!} \ll A^K 4^\nu \nu! \sum_{\mu=0}^\infty \frac{(AK)^\mu}{\mu!} \sum_{k=0}^\infty \frac{(AK)^k}{k!} \ll (AK)^K,
\end{split}
\]
and \eqref{L3-19/3} follows.

\smallskip
In order to compute the degree of $Q_\nu(s)$ we first recall that $H_F(0)=d_F=2$. From the definition of $R_\nu(s)$, see \eqref{3-0bis}, and recalling that the leading term of $B_n(x)$ is $x^n$, see (3) of Sect.1.13 of Bateman's Project \cite{EMOT/1953}, we have that the leading term of $R_\nu(s)$ is
\[
s^{\nu+1}\{(-2)^{\nu+1} +(-1)^\nu - (-1)^{\nu+1}\} = s^{\nu+1}\{(-2)^{\nu+1} +2(-1)^\nu\},
\]
hence $\deg R_\nu(s)=\nu+1$ for $\nu\geq 1$. From the definition of $V_\nu(s)$, see \eqref{3-1}, and the previous assertion we easily see that the degree of $V_\nu(s)$ is given by the single term on the right hand side of \eqref{3-1} arising for $m=\mu$, thus $\deg V_\nu(s)=2\nu$. To prove that $\deg Q_\nu(s)=2\nu$ we first note that $\deg A_{\mu,\nu}(s) = \nu-\mu$ and $C_{\nu,\nu}\neq0$, thanks to Lemmas 3.14 and 3.12. Hence the degree of $Q_\nu(s)$ is given by the term with $\mu=\nu$ in \eqref{L3-15/0}, and the lemma follows. \fine

\bigskip
{\bf 3. Proof of Theorem 2.} We follow the notation in subsection 1. Moreover, let $z_X=\frac{1}{X}+2\pi i\alpha$ with a large $X>0$. At the beginning we keep open the value of the sufficiently large constants $c_0,c_1,c_2,N$, and we add conditions on them when required. Writing
\[
F_X(s,\alpha) = \sum_{n=1}^\infty\frac{a(n)}{n^s}\exp(-nz_X),
\]
for $\sigma<2$ we have
\begin{equation}
\label{P1}
F_X(s,\alpha) = \frac{1}{2\pi i}\int_{(2-\sigma)} F(s+w)\Gamma(w)z_X^{-w}\d w = \frac{1}{2\pi i}\int_{\LL(s)} F(s+w)\Gamma(w)z_X^{-w}\d w,
\end{equation}
since the poles of $F(s+w)\Gamma(s)$ lie to the left of the path $\LL(s)$. If $w\in\LL_{-\infty}(s)$ then $\Re(s+w)\geq -c_0$, hence 
$F(s+w)\ll |s+w|^c$ for some $c>0$ since $F(s)$ has polynomial growth on vertical strips. If in addition $v<t_0$, then $\Re(s+w)=c_0>1$ and hence $F(s+w)\ll 1$. Moreover, still for $w\in\LL_{-\infty}(s)$, we have
\[
|z_X^{-w}|=|z_X|^{-u}\exp(v(\pi/2+O(1/X)))
\]
\[
|\Gamma(w)| \leq \Gamma(|\sigma|+c_0) \ll (|\sigma|+c_0)^{|\sigma|+c_0} \ll A^{|\sigma|}|\sigma|^{|\sigma|}+1.
\]
Therefore the contribution of $\LL_{-\infty}(s)$ to \eqref{P1} is
\[
\frac{1}{2\pi i}\int_{\LL_{-\infty}(s)} F(s+w)\Gamma(w)z_X^{-w}\d w  \ll A^{|s|}|\sigma|^{|\sigma|}+1
\]
for some $A>0$. As a consequence, for $\sigma<2$ and any fixed $\alpha>0$
\begin{equation}
\label{P2}
\begin{split}
F_X(s,\alpha) &= \frac{1}{2\pi i}\int_{\LL_{\infty}(s)} F(s+w)\Gamma(w)z_X^{-w}\d w + O(A^{|s|}|\sigma|^{|\sigma|}+1) \\
&= \I_X(s,\alpha) + O(A^{|s|}|\sigma|^{|\sigma|}+1),
\end{split}
\end{equation}
say, uniformly as $X\to\infty$. Note in passing that $\I_X(s,\alpha)$ is not holomorphic in $s$, since the path $\LL_{\infty}(s)$ starts at $-\sigma-c_0+it_0$, which is not holomorphic in $s$.

\medskip
In order to study the integral $\I_X(s,\alpha)$ we apply the functional equation of $F(s)$ and the reflection formula of $\Gamma(s)$, thus getting
\[
\I_X(s,\alpha) = \omega Q^{1-2s} \frac{1}{2\pi i}\int_{\LL_{\infty}(s)} \bar{F}(1-s-w)G(s,w)S(s,w)(Q^2z_X)^{-w}\d w.
\]
Replacing $S(s,w)$ by $-ie(-\xi_F/4)e^{-\pi is}$ and estimating the error by Lemma 3.2 and then by Lemma 3.1 we obtain
\begin{equation}
\label{P3}
\begin{split}
\I_X(s,\alpha) &= -i\omega e(-\xi_F/4)Q^{1-2s}e^{-\pi is} \frac{1}{2\pi i}\int_{\LL_{\infty}(s)} \bar{F}(1-s-w)G(s,w)(Q^2z_X)^{-w}\d w \\
& \hskip1cm + O(A^{|s|} \int_{t_0}^\infty|G(s,-\sigma-c_0+iv)| |(Q^2z_X)^{\sigma+c_0-iv}|e^{-cv}\d v) \\
&= \J_X(s,\alpha) + O(A^{|s|}|\sigma|^{|\sigma|}+1),
\end{split}
\end{equation}
say, uniformly as $X\to\infty$.

\medskip
Now, roughly speaking, we reduce $G(s,w)$ in \eqref{P3} to a single $\Gamma$-factor by means of the uniform version of the Stirling formula in \cite{Ka-Pe/USF}. Clearly
\[
\begin{split}
\log G(s,w) &= (1-r)\log 2\pi -\log \Gamma(-w+1) \\
&+\sum_{j=1}^r\big\{ \log\Gamma(-\lambda_jw+\lambda_j(1-s)+\bar{\mu}_j) +\log\Gamma(-\lambda_jw+1-\lambda_js-\mu_j)\big\},
\end{split}
\]
and we apply the Theorem in \cite{Ka-Pe/USF} with the choices
\[
(z,s) = (-w,1), \quad (z,s) = (-\lambda_jw,\lambda_j(1-s)+\bar{\mu}_j), \quad (z,s) = (-\lambda_jw,1-\lambda_js-\mu_j)
\]
to the above three $\log\Gamma$-terms, respectively. Since $w\in\LL_{\infty}(s)$, it is easy to verify that the hypotheses of the above quoted Theorem are satisfied provided the constants in the definition of $\LL_{\infty}(s)$ are large enough. Hence for $N\leq |s|+c$ we get
\[
\begin{split}
\log G(s,w) &= (1-r)\log 2\pi + \sum_{j=1}^r\big\{(\lambda_j(1-s-w)+\bar{\mu}_j-\frac12)\log(-w) \\
&\qquad+ (\lambda_j(1-s-w)+\bar{\mu}_j-\frac12)\log\lambda_j + (\frac12-\lambda_j(s+w)-\mu_j)\log(-w) \\
&\qquad+ (\frac12-\lambda_j(s+w)-\mu_j)\log\lambda_j + 2\lambda_jw +\log2\pi \big \} \\
&\qquad+ \sum_{j=1}^r\sum_{\nu=1}^N \frac{(-1)^{\nu+1}}{\nu(\nu+1)} \big(B_{\nu+1}(\lambda_j(1-s)+\bar{\mu}_j) + B_{\nu+1}(1-\lambda_js-\mu_j)\big) \big(-\frac{1}{\lambda_jw}\big)^\nu \\
&\qquad- (\frac12-w)\log(-w) -w -\frac12\log2\pi -\sum_{\nu=1}^N \frac{(-1)^{\nu+1}}{\nu(\nu+1)} B_{\nu+1}(1) \big(-\frac{1}{w}\big)^\nu \\
&\qquad+O\big(\frac{(c(|s|+1))^{N+2}}{|w|^{N+1}}\big) \\
&= (\frac12-2s-w-i\theta_F)\log(-w) + w + \frac12\log2\pi +(1-2s-2w)\log\prod_{j=1}^r\lambda_j^{\lambda_j} \\
&\qquad+ \log\prod_{j=1}^r\lambda_j^{-2i\Im\mu_j} + \sum_{\nu=1}^N \frac{P_\nu(s)}{\nu(\nu+1)} \frac{1}{w^\nu} + O\big(\frac{(c(|s|+1))^{N+2}}{|w|^{N+1}}\big),
\end{split}
\]
say, where by Lemma 3.3
\[
\begin{split}
P_\nu(s) &= B_{\nu+1}(1) - \sum_{j=1}^r \frac{B_{\nu+1}(\lambda_j(1-s)+\bar{\mu}_j) + B_{\nu+1}(1-\lambda_js-\mu_j)}{\lambda_j^\nu} \\
&= R_\nu(s) - B_{\nu+1}(1-2s-i\theta_F).
\end{split}
\]
On the other hand, again from the Theorem in \cite{Ka-Pe/USF} but with the choice $(z,s) = (-w,1-2s-i\theta_F)$, we obtain
\[
\begin{split}
\log\Gamma(1-2s-w-i\theta_F) &= (\frac12 -2s-w-i\theta_F)\log(-w) + w + \frac12\log2\pi \\
&- \sum_{\nu=1}^N \frac{B_{\nu+1}(1-2s-i\theta_F)}{\nu(\nu+1)} \frac{1}{w^\nu} + O\big(\frac{(c(|s|+1))^{N+2}}{|w|^{N+1}}\big).
\end{split} 
\]
Therefore, writing $\beta= \prod_{j=1}^r\lambda_j^{2\lambda_j}$ we have
\begin{equation}
\label{P3bis}
\begin{split}
\log G(s,w) &= \log\Gamma(1-2s-w-i\theta_F) + (\frac12-s-w)\log\beta + \log\prod_{j=1}^r\lambda_j^{-2i\Im\mu_j} \\
&+ \sum_{\nu=1}^N \frac{R_\nu(s)}{\nu(\nu+1)} \frac{1}{w^\nu} + O\big(\frac{(c(|s|+1))^{N+2}}{|w|^{N+1}}\big).
\end{split} 
\end{equation}
Since $O\big(\frac{(c(|s|+1))^{N+1}}{|w|^{N+1}}\big) < 1$ for $w\in\LL_\infty(s)$ provided the constants in the definition of $\LL_\infty(s)$ are sufficiently large, we have
\[
e^{O\big(\frac{(c(|s|+1))^{N+2}}{|w|^{N+1}}\big)} = 1 + \sum_{k=1}^\infty \frac{(c(|s|+1))^k}{k!} \left(O\big(\frac{(c(|s|+1))^{N+1}}{|w|^{N+1}}\big)\right)^k = 1 + O\big(A^{|s|}\frac{(c(|s|+1))^{N+1}}{|w|^{N+1}}\big),
\]
hence from \eqref{P3bis} we get
\begin{equation}
\label{P4}
\begin{split}
G(s,w) &= \Gamma(1-2s-w-i\theta_F) \beta^{\frac12-s-w} \prod_{j=1}^r\lambda_j^{-2i\Im\mu_j} \exp\big(\sum_{\nu=1}^N \frac{R_\nu(s)}{\nu(\nu+1)} \frac{1}{w^\nu}\big) \\
&\times \big(1+ O\big(A^{|s|}\frac{(c(|s|+1))^{N+1}}{|w|^{N+1}}\big)\big).
\end{split} 
\end{equation}

\medskip
Replacing in $\J_X(s,\alpha)$ (see \eqref{P3}) $G(s,w)$ by its main term obtained in \eqref{P4} causes an error of size
\[
\ll A^{|s|} \int_{\LL_\infty(s)} |\Gamma(1-2s-w-i\theta_F)| \exp\big(\Re\big(\sum_{\nu=1}^N \frac{R_\nu(s)}{\nu(\nu+1)} \frac{1}{w^\nu}\big)\big) \big|(Q^2z_X)^{-w}\big| \frac{(c(|s|+1))^{N+3}}{|w|^{N+1}} |\d w|,
\]
and by Lemmas 3.4 and 3.9 this is
\[
\ll A^{|s|}\big(c(|s|+1)\big)^{N+3} \int_{t_0}^\infty e^{-\pi v/2} e^{\pi v/2} v^{-\sigma+2c_0-N} \d v
\]
uniformly in $X$. Hence if $N\geq -\sigma+3c_0$ the integral converges and is $\ll 1$. Moreover, if $N\leq -\sigma+c$ (with any fixed $c>3c_0$) then $(2|s|)^{N+3}+1 \ll A^{|s|}(|s|^{|\sigma|}+1)$. Therefore, from \eqref{P3}, \eqref{P4} and recalling the definition of conductor $q_F$ and root number $\omega^*_F$, for such $N$'s we have
\[
\begin{split}
\J_X(s,\alpha) &= \omega^*_F \big(\frac{q_F}{4\pi^2}\big)^{1/2 -s -i\theta_F/2} e(-\frac12(s+i\frac{\theta_F}{2}) \frac{1}{2\pi i} \int_{\LL_{\infty}(s)} \bar{F}(1-s-w) \ \times \\
&\qquad \times \Gamma(1-2s-w-i\theta_F) \exp\big(\sum_{\nu=1}^N \frac{R_\nu(s)}{\nu(\nu+1)} \frac{1}{w^\nu}\big)\big(\frac{q_Fz_X}{4\pi^2}\big)^{-w} \d w + O\big(A^{|s|} (|s|^{|\sigma|}+1)\big).
\end{split}
\]
Hence, by the substitution $1-2s-w-i\theta_F \to w$ in the above integral, for
\begin{equation}
\label{P5}
-\sigma+3c_0 \leq N \leq -\sigma+c
\end{equation}
we obtain
\begin{equation}
\label{P6}
\begin{split}
\J_X(s,\alpha) &= -i\omega^*_F \big(\sqrt{q_F}\alpha - i\frac{\sqrt{q_F}}{2\pi X}\big)^{2s-1+i\theta_F} \frac{1}{2\pi i} \int_{\LL^*_{-\infty}(s)} \bar{F}(s+w+i\theta_F) \Gamma(w) \ \times \\
&\qquad \times \exp\big(\sum_{\nu=1}^N \frac{(-1)^\nu R_\nu(s)}{\nu(\nu+1)} \frac{1}{(w+2s-1+i\theta_F)^\nu}\big)\big(\frac{q_Fz_X}{4\pi^2}\big)^w \d w + O\big(A^{|s|} (|s|^{|\sigma|}+1)\big),
\end{split}
\end{equation}
uniformly as $X\to\infty$.

\medskip
Now we use Lemma 3.15 to replace the term $\exp(\sum_{\nu=1}^N...)$ in the above integral by the sum involving the polynomials $Q_\nu(s)$ with $\nu=0,...,N$. Hence we need that $c$ in \eqref{P5} satisfies $c\leq c_2/3$. Since $|F(s+w+i\theta_F)|\ll 1$ for $w\in\LL^*_{-\infty}(s)$, this causes a further error of size
\[
\begin{split}
&\ll A^{|s|}(|s|^{|\sigma|}+1) \int_{\LL^*_{-\infty}(s)} \frac{|\Gamma(w)|}{|w(w-1)\cdots(w-N)|} e^{\pi|v|/2} |\d w| \\
&\ll A^{|s|}(|s|^{|\sigma|}+1) \int_{\LL^*_{-\infty}(s)} \frac{|\Gamma(w-N)|e^{\pi|v|/2}}{|w|} |\d w|.
\end{split}
\]
Moreover, we also want $N$ such that $\Re(w-N)\leq 0$ for $w\in\LL^*_{-\infty}(s)$, i.e. we choose
\begin{equation}
\label{P6bis}
N=[-\sigma]+k
\end{equation}
with a sufficiently large positive integer $k$ satisfying \eqref{P5} with $c\leq c_2/3$; this can be done by suitably choosing $c_0$ and $c_2$. With such a choice of $N$ we have $|\Gamma(w-N)|\ll e^{-\pi|v|/2}/|v|^{1/2}$, hence the integral is $\ll 1$. Therefore \eqref{P6} becomes, uniformly as $X\to\infty$,
\begin{equation}
\label{P7}
\begin{split}
\J_X(s,\alpha) &= -i\omega^*_F \big(\sqrt{q_F}\alpha - i\frac{\sqrt{q_F}}{2\pi X}\big)^{2s-1+i\theta_F} \sum_{\nu=0}^N Q_\nu(s) \frac{1}{2\pi i} \int_{\LL^*_{-\infty}(s)} \bar{F}(s+w+i\theta_F)  \ \times \\
&\qquad \times \Gamma(w-\nu) \big(\frac{q_Fz_X}{4\pi^2}\big)^w \d w + O\big(A^{|s|} (|s|^{|\sigma|}+1)\big).
\end{split}
\end{equation}

\medskip
Replacing the path of integration in \eqref{P7} by the whole path $\LL^*(s)$ causes an error which, since $\nu\leq N$, by Lemma 3.18 is of size
\[
\ll A^{|s|} \sum_{\nu=0}^N \frac{(|s|+1)^{2\nu}}{\nu!} \int_{\LL_\infty^*(s)} |F(s+w+i\theta_F) \Gamma(w-\nu) \big(\frac{q_Fz_X}{4\pi^2}\big)^w| |\d w|.
\]
For $w\in\LL_\infty(s)$ we have $\Re(s+w)=O(1)$ and hence $F(s+w+i\theta_F)\ll (|s|+|w|+1)^A$. Moreover, for $0\leq \nu\leq N$ we have $-c \leq \Re(w-\nu) \leq |\sigma|-\nu+c'$ and hence
\[
\Gamma(w-\nu) \ll (c(|\sigma|+1))^{|\sigma|-\nu+c''}.
\]
Further, $|z_X^w| \ll A^{|\sigma|} e^{-v\arg z_X}$, thus the above mentioned error is
\[
\ll A^{|s|} \sum_{\nu=0}^N \frac{(|s|+1)^{|\sigma|+\nu}}{\nu!} \int_{-t_0^*(s)}^\infty (|v|+1)^A e^{-v\arg{z_X}}\d v \ll A^{|s|}(|s|+1)^{|\sigma|}.
\]
We also have by Cauchy's theorem that for $0\leq \nu\leq N$
\[
\begin{split}
&\frac{1}{2\pi i} \int_{\LL^*(s)} \bar{F}(s+w+i\theta_F) \Gamma(w-\nu) \big(\frac{q_Fz_X}{4\pi^2}\big)^w \d w \\
&\hskip-.5cm = \frac{1}{2\pi i} \int_{|\sigma|+\nu+2-i\infty}^{|\sigma|+\nu+2+i\infty}  \bar{F}(s+w+i\theta_F) \Gamma(w-\nu) \big(\frac{q_Fz_X}{4\pi^2}\big)^w \d w \\
&\hskip-.5cm = \sum_{n=1}^\infty \frac{\overline{a(n)}}{n^{s+i\theta_F}} \big(\frac{q_Fz_X}{4\pi^2n}\big)^\nu \frac{1}{2\pi i} \int_{|\sigma|+\nu+2-i\infty}^{|\sigma|+\nu+2+i\infty}  \Gamma(w-\nu) \big(\frac{q_Fz_X}{4\pi^2n}\big)^{w-\nu} \d w \\
&\hskip-.5cm =  \big(\frac{q_F}{4\pi^2X}+i\frac{q_F\alpha}{2\pi}\big)^\nu \sum_{n=1}^\infty \frac{\overline{a(n)}}{n^{s+\nu+i\theta_F}} \exp(-\frac{4\pi^2}{q_Fz_X}n).
\end{split}
\]
Consequently, \eqref{P7} becomes
\begin{equation}
\label{P8}
\begin{split}
\J_X(s,\alpha) &= -i\omega^*_F \big(\sqrt{q_F}\alpha - i\frac{\sqrt{q_F}}{2\pi X}\big)^{2s-1+i\theta_F} \sum_{\nu=0}^N 
 \big(\frac{q_F}{4\pi^2X}+i\frac{q_F\alpha}{2\pi}\big)^\nu Q_\nu(s) \times \\
&\hskip.8cm \times \sum_{n=1}^\infty \frac{\overline{a(n)}}{n^{s+\nu+i\theta_F}} \exp(-\frac{4\pi^2}{q_Fz_X}n) +O(A^{|s|}(|s|+1)^{|\sigma|})
\end{split}
\end{equation}
uniformly as $X\to\infty$. Since
\[
-\frac{4\pi^2}{q_Fz_X} = \frac{2\pi i}{q_F\alpha} - \frac{1}{q_F\alpha^2}\frac{1}{X+O(1)},
\]
the series in \eqref{P8} is absolutely convergent for all $s$, for every $\nu$.

\medskip
The next step is to make the range of summation of $\nu$ in \eqref{P8} independent of $\sigma$ (recall that $N$ depends on $\sigma$, see \eqref{P6bis}). Let $K>0$ be a large integer and $|s|< 2K$, $\sigma> -K+1/2$. Depending on the relative sizes of $N$ and $K$, we add to or withdrow from \eqref{P8} the terms with $\nu$ between $N+1$ and $K$ or between $K+1$ and $N$, respectively. In both cases we have that $\sigma+\nu> 3/2$ for such $\nu$'s (call them $\nu\in\X$), hence from Lemma 3.19 we deduce that
\begin{equation}
\label{P9}
\begin{split}
-i\omega^*_F &\big(\sqrt{q_F}\alpha - i\frac{\sqrt{q_F}}{2\pi X}\big)^{2s-1+i\theta_F} \sum_{\nu\in\X} \big(\frac{q_F}{4\pi^2X}+i\frac{q_F\alpha}{2\pi}\big)^\nu Q_\nu(s) \times \\
&\times \sum_{n=1}^\infty \frac{\overline{a(n)}}{n^{s+\nu+i\theta_F}} \exp(-\frac{4\pi^2}{q_Fz_X}n) \ll A^K \sum_{\nu\in\X} A^\nu (AK)^K \ll (A'K)^K
\end{split}
\end{equation}
uniformly in $X$. From \eqref{P2}, \eqref{P3}, \eqref{P8} and \eqref{P9} we therefore obtain that for $-K+1/2<\sigma<2$ and $|s|< 2K$
\begin{equation}
\label{P10}
\begin{split}
F_X(s,\alpha) = &-i\omega^*_F \big(\sqrt{q_F}\alpha - i\frac{\sqrt{q_F}}{2\pi X}\big)^{2s-1+i\theta_F} \sum_{\nu=0}^K \big(\frac{q_F}{4\pi^2X}+i\frac{q_F\alpha}{2\pi}\big)^\nu \times \\
&\times Q_\nu(s) F^*_X(s+\nu+i\theta_F,\alpha) + H_X(s,\alpha),
\end{split}
\end{equation}
where
\[
F^*_X(s,\alpha) = \sum_{n=1}^\infty \frac{\overline{a(n)}}{n^{s+\nu+i\theta_F}} \exp(-\frac{4\pi^2}{q_Fz_X}n)
\]
and 
\[
H_X(s,\alpha) \ll (AK)^K
\]
uniformly as $X\to\infty$. Moreover, since $F_X(s,\alpha)$, $F^*_X(s,\alpha)$ and $Q_\nu(s)$ are entire functions, $H_X(s,\alpha)$ is also entire. Further, from \eqref{P10} we have that for $1<\sigma<2$
\[
\lim_{X\to\infty} H_X(s,\alpha) = H(s,\alpha)
\]
exists and is holomorphic since this is clearly true for $F_X(s,\alpha)$ and $F^*_X(s,\alpha)$. For $1<\sigma<2$
we also have that 
\begin{equation}
\label{P11}
H(s,\alpha) = F(s,\alpha) + i\omega^*_F \big(\sqrt{q_F}\alpha\big)^{2s-1+i\theta_F} \sum_{\nu=0}^K \big(i\frac{q_F\alpha}{2\pi}\big)^\nu Q_\nu(s) \bar{F}(s+\nu+i\theta_F,-\frac{1}{q_F\alpha}).
\end{equation}
Hence by Vitali's convergence theorem, see Section 5.21 of Titchmarsh \cite{Tit/1939}, the limit function $H(s,\alpha)$ exists and is holomorphic for $-K+1/2<\sigma<2$ and $|s|<2K$, and satisfies
\[
H(s,\alpha) \ll (AK)^K.
\]
This provides analytic continuation and bounds for the right hand side of \eqref{P11}. Therefore Theorem 2 follows, recalling Lemmas 3.18 and 3.19.

\bigskip
\section{Proof of Theorem 3 and Corollary}  

\smallskip
Let $m_F$ denote the order of pole of $F(s)$ at $s=1$ and let $Q_\nu(s)$ be the polynomials in Theorem 2. We have

\medskip
{\bf Lemma 4.1.} {\sl Let $F\in\S_2^\sharp$ with $q_F=1$. For $\nu\geq 1$ we have that $(s+\nu-1)^{m_F}$ divides $Q_\nu(s)$. Moreover, $\theta_F=0$ if $m_F>0$.}

\medskip
{\sl Proof.} We may clearly assume that $m_F>0$. From Theorem 2 with $\alpha=1$ we get
\begin{equation}
\label{4-0}
F(s) = -i\omega^*_F \bar{F}(s+i\theta_F) + H(s),
\end{equation}
where $H(s)$ is holomorphic for $\sigma>1/2$. Hence $\theta_F=0$ since $m_F>0$. Again from Theorem 2 with $\alpha=1$ we deduce that
\begin{equation}
\label{4-1}
F(s) = -i\omega^*_F \sum_{\nu=0}^K \big(\frac{i}{2\pi}\big)^\nu Q_\nu(s)\bar{F}(s+\nu) + H_K(s),
\end{equation}
where $K>0$ is an arbitrarily large integer and $H_K(s)=H_K(s,1)$ is as in Theorem 2. Given $1\leq \nu_0\leq K$, it is clear that all terms in \eqref{4-1} are holomorphic at $s=1-\nu_0$, except possibly for $Q_{\nu_0}(s)\bar{F}(s+\nu_0)$. Therefore this term must be holomorphic as well, hence $Q_{\nu_0}(s)$ has a zero of order at least $m_F$ at $s=1-\nu_0$, and the result follows. \fine

\medskip
{\bf Lemma 4.2.} {\sl Let $F\in\S_2^\sharp$ with $q_F=1$. For $q\geq 1$ and $1\leq a\leq q$ with $(a,q)=1$ the function $(s-1)^{m_F}F(s,a/q)$ is entire.}

\medskip
{\sl Proof.} In view of Lemma 4.1, for $\nu\geq 1$ we define the polynomials $P_\nu(s)$ by
\[
Q_\nu(s) = (s+\nu-1)^{m_F} P_\nu(s). 
\]
For $q=1$ the result is obvious, and we proceed by induction on $q$. Assume the result true up to $q-1$ and apply Theorem 2 with $\alpha=a/q$, $1\leq a<q$ and $(a,q)=1$, thus getting
\[
\begin{split}
(s-1)^{m_F}F(s,a/q) &= -i\omega_F^*  \big(\frac{a}{q}\big)^{2s-1+i\theta_F} \big\{(s-1)^{m_F}\bar{F}(s+i\theta_F,-q/a) + (s-1)^{m_F} \\
&\times \sum_{\nu=1}^K \big(\frac{ia}{2\pi q}\big)^\nu P_\nu(s) (s+\nu-1)^{m_F} \bar{F}(s+\nu+i\theta_F,-q/a) + H_K(s,a/q).
\end{split}
\]
Since $a\leq q-1$, by the inductive hypothesis all terms on the right hand side are holomorphic in the domain where $H_K(s,a/q)$ is holomorphic (remember that $\theta_F=0$ if $m_F>0$ by Lemma 4.1). The result follows since $K$ is arbitrarily large. \fine

\medskip
In order to prove Theorem 3 we note that thanks to Lemma 4.2 we only need to show that $F(s,a/q)$ is of the proper size as $\sigma\to-\infty$. We proceed by induction on $q$ and observe that for $q=1$ the result follows from Lemma 2.1 since $F\in M(2,1)$. Assume now the result true up to $q-1$ and apply Theorem 2 with $\alpha=a/q$, thus getting for arbitrarily fixed $K>0$, $A,B\in\RR$ and a suitable $C>0$ (whose value will not necessary be the same at each occurrence) that
\[
F(s,a/q) = -i\omega_F^* \big(\frac{a}{q}\big)^{2s-1+i\theta_F} \sum_{\nu=0}^K \big(\frac{ia}{2\pi q}\big)^\nu Q_\nu(s) \bar{F}(s+\nu+i\theta_F,-q/a) + O(C^KK^K)
\]
for $\sigma>-K+1/2$ and $|s|<2K$, uniformly for $A\leq t\leq B$. Choosing $K=[|\sigma|]+2$ and letting $\sigma\to-\infty$ we obtain
\[
F(s,a/q) \ll \big(\frac{a}{q}\big)^{2|\sigma|} \sum_{\nu=0}^K \big(\frac{a}{2\pi q}\big)^\nu |Q_\nu(s)| |F(s+\nu+i\theta_F,-q/a)| + C^{|\sigma|} |\sigma|^{|\sigma|}.
\]
Hence by Lemma 3.18 and the inductive hypothesis we have
\[
\begin{split}
F(s,a/q) &\ll \big(\frac{q}{a}\big)^{2|\sigma|} \sum_{\nu=0}^{[|\sigma|]+2} \big(\frac{a}{2\pi q}\big)^\nu \frac{(C|\sigma|)^{2\nu}}{\nu!} (|\sigma|+2-\nu)^{2(|\sigma|-\nu)} \big(\frac{a}{2\pi e}\big)^{2(|\sigma|-\nu)} |\sigma|^C+ C^{|\sigma|} |\sigma|^{|\sigma|} \\
&\ll |\sigma|^{2|\sigma|} \big(\frac{q}{2\pi e}\big)^{2|\sigma|} |\sigma|^C \sum_{\nu=0}^{[|\sigma|]+2} \frac{C^\nu}{\nu!} + C^{|\sigma|} |\sigma|^{|\sigma|} \ll |\sigma|^{2|\sigma|} \big(\frac{q}{2\pi e}\big)^{2|\sigma|} |\sigma|^C.
\end{split}
\]
Therefore $F(s,a/q)$ belongs to $M(2,q^2)$, and Theorem 3 follows. \fine

\bigskip
In view of Theorems 1 and 3, to prove the Corollary we need to show that every $F\in\S_2$ satisfies $N_F(\sigma,T) = o(T)$ for any fixed $\sigma>1/2$. Actually, standard techniques (see Ch.12 of Montgomery \cite{Mon/1971}) allow to show the following sharper result: if $F\in\S_2$ then for every $\epsilon>0$ and every fixed $\sigma>1/2$
\[
N_F(\sigma,T) \ll T^{3/2-\sigma+\epsilon}.
\]
We only outline the main points in the proof. Let $a(n)$ be the Dirichlet coefficients of $F(s)$ and let $\mu_F(n)$ denote its inverse. By Lemma 1 of \cite{Ka-Pe/2003} we have that for every $\epsilon>0$ there exists an integer $M=M(\epsilon)$ such that $\mu_F(n) \ll n^\epsilon$ for $(n,M)=1$. Moreover, for $\sigma>1/2$ the functions $F(s)$ and 
\[
F_M(s) = F(s) \prod_{p|M}F_p(s)^{-1} = \sum_{(n,M)=1} \frac{a(n)}{n^s}
\]
have the same zeros. Writing $L=\log T$ and
\[
G(s) = F_M(s) \sum_{\substack{n\leq TL^2 \\ (n,M)=1}} \frac{\mu_F(n)}{n^s} = F_M(s)M(s),
\]
say, we have
\[
G(s) =1 + \sum_{\substack{n>TL^2 \\ (n,M)=1}}\frac{c(n)}{n^s}, \hskip2cm c(n) = \sum_{\substack{d|n \\ d\leq TL^2}} \mu_F(d)a(n/d) \ll n^\epsilon.
\]
Now we apply Montgomery's zero detecting method. First, for every zero $\rho=\beta+i\gamma$ of $F(s)$ with $\beta>\sigma>1/2$ we obtain
\begin{equation}
\label{4-2}
I_\rho = \frac{1}{2\pi i} \int_{2-i\infty}^{2+i\infty} G(\rho+w)\Gamma(w)T^w \d w = 1 + O(\frac{1}{T}).
\end{equation}
On the other hand, shifting the line of integration to $\Re w=1/2-\beta$ we get
\begin{equation}
\label{4-3}
I_\rho \ll T^{1/2-\sigma} \int_{\gamma-cL}^{\gamma+cL} |F(\frac12 ++it) M(\frac12 ++it)| \d t +  O(\frac{1}{T})
\end{equation}
with a suitable $c>0$. Summing over representatives of zeros in small rectangles, from \eqref{4-2}, \eqref{4-3}, the Cauchy-Schwarz inequality and the mean-value theorem for Dirichlet polynomials we obtain
\[
N_F(\sigma,T) \ll T^{1/2-\sigma}L^3 \big(\int_{T/2}^{2T} |F(\frac12 ++it)|^2 \d t\big)^{1/2} \big(T\sum_{\substack{n\leq TL^2 \\ (n,M)=1}}\frac{|\mu_F(n)|^2}{n}\big)^{1/2} \ll T^{3/2-\sigma+\epsilon}.
\]
Here we used the bound
\[
\int_{T/2}^{2T} |F(\frac12 ++it)|^2 \d t \ll T^{1+\epsilon},
\]
which follows by standard arguments from the approximate functional equation in Chan- drasekharan-Narasimhan \cite{Ch-Na/1963} for $L$-functions of degree $d=2$ ($A=1$ in the notation of \cite{Ch-Na/1963}). The Corollary is therefore proved.

\bigskip
\section{Proof of Theorem 4}

\smallskip
We need further notation. For $F\in\S^\sharp_2$ we write
\[
\alpha_F(a/q) = \lim_{s\to 1} (s-1)^{m_F} F(s,a/q). 
\]
Moreover, let $\alpha_F=\alpha_F(1)$ and $\lambda_F=-i\omega_F^*$. Note that $\alpha_F\neq 0$ if $m_F>0$.

\medskip
{\bf Lemma 5.1.} {\sl Let $F\in\S_2^\sharp$ with $q_F=1$ and $m_F>0$. Then}
\[
\text{{\sl i)}} \  \alpha_{\bar{F}}=\overline{\alpha_F},\hskip1cm \text{{\sl ii)}} \ m_F\leq 2, \hskip1cm \text{{\sl iii)}} \ \alpha_F=\lambda_F \overline{\alpha_F}, \hskip1cm \text{{\sl iv)}} \ \lambda_{\bar{F}} = \overline{\lambda_F}. 
\]

\medskip
{\sl Proof.} {\sl i)} is trivial. By Lemmas 3.18 and 4.1 we have $m_F\leq \deg Q_1(s)\leq 2$, hence {\sl ii)} follows. Multiplying both sides of \eqref{4-0} by $(s-1)^{m_F}$ and letting $s\to 1$ we obtain 
\[
\alpha_F = \lambda_F\alpha_{\bar{F}}= \lambda_F\overline{\alpha_F},
\]
and {\sl iii)} follows. Finally, applying {\sl iii)} to $\bar{F}(s)$, thanks to {\sl i)}  we get
\[
\alpha_F = \overline{\alpha_{\bar{F}}}  = \overline{\lambda_{\bar{F}}} \overline{\alpha_F},
\]
and {\sl iv)} follows comparing with {\sl iii)}. \fine

\medskip
{\bf Lemma 5.2.} {\sl Let $F\in\S_2^\sharp$ with $q_F=1$ and $m_F>0$. Then for $q\geq 1$ and $1\leq a\leq q$ with $(a,q)=1$ we have}
\[
\alpha_F(a/q) = \frac{\alpha_F}{q}.
\]

\medskip
{\sl Proof.} We proceed by induction, the case $q=1$ being trivial. Recalling that $\theta_F=0$ in this case, from Theorem 2 with $\alpha=a/q$ we get
\[
F(s,a/q) = \lambda_F\big(\frac{a}{q}\big)^{2s-1} \overline{F(\bar{s},q/a)}+ H(s)
\]
with $H(s)$ holomorphic for $\sigma>1/2$. Multiplying both sides by $(s-1)^{m_F}$ and letting $s\to 1$ we obtain, thanks to the inductive hypothesis, that
\[
\alpha_F(a/q) = \lambda_F \frac{a}{q} \overline{\alpha_F(q/a)} = \lambda_F \frac{a}{q} \frac{\overline{\alpha_F}}{a} = \frac{\lambda_F\overline{\alpha_F}}{q}.
\]
The result follows now by {\sl iii)} of Lemma 5.1. \fine

\medskip
When $m_F=2$ we write
\[
\begin{split}
F(s) &= \frac{\alpha_F}{(s-1)^2} + \frac{\alpha_F\beta_F}{s-1} + \dots \\
F(s,a/q) &= \frac{\alpha_F(a/q)}{(s-1)^2} + \frac{\alpha_F(a/q)\beta_F(a/q)}{s-1} + \dots
\end{split}
\]

\medskip
{\bf Lemma 5.3.} {\sl Let $F\in\S_2^\sharp$ with $q_F=1$ and $m_F=2$. Then $\beta_F\in\RR$, and for $q\geq 1$ and $1\leq a\leq q$ with $(a,q)=1$ we have} 
\[
\beta_F(a/q) = \beta_F -2\log q.
\]

\medskip
{\sl Proof.} In order to prove that $\beta_F\in\RR$ we start again with Theorem 2 with $\alpha=1$ (see \eqref{4-0}), hence
\[
F(s) = \lambda_F\bar{F}(s) + H(s)
\]
with $H(s)$ holomorphic for $\sigma>1/2$. This gives
\[
 \frac{\alpha_F}{(s-1)^2} + \frac{\alpha_F\beta_F}{s-1} + \dots =  \frac{\lambda_F\alpha_{\bar{F}}}{(s-1)^2} + \frac{\lambda_F\alpha_{\bar{F}}\beta_{\bar{F}}}{s-1} + \dots,
\]
hence thanks to Lemma 5.1 we deduce that $\beta_F=\beta_{\bar{F}}$. But for $s\in\RR$ we have
\[
\beta_{\bar{F}} = \lim_{s\to 1} \big(\frac{(s-1)\bar{F}(s)}{\alpha_{\bar{F}}}-\frac{1}{s-1}\big) = \lim_{s\to 1} \overline{\big(\frac{(s-1)F(s)}{\alpha_{F}}-\frac{1}{s-1}\big)} = \overline{\beta_F},
\]
and the first assertion follows. Now we prove the second assertion by induction on $q$, the case $q=1$ being trivial. Once again from Theorem 2 with $\alpha = a/q$, writing the first terms of the Laurent expansion at $s=1$ of the right hand side, we get
\[
F(s,a/q) = \lambda_F(\frac{a}{q} + 2\frac{a}{q}\log\frac{a}{q} (s-1) + \dots) (\frac{\overline{\alpha_F}}{a(s-1)^2} + \frac{\overline{\alpha_F}\overline{\beta_F(q/a)}}{a(s-1)} + \dots) + H(s)
\]
with $H(s)$ holomorphic for $\sigma>1/2$. By the inductive assumption we have $\overline{\beta_F(q/a)} = \beta_F -2\log a$. Hence from the previous equation, thanks to Lemmas 5.1 and 5.2 we obtain
\[
\begin{split}
F(s,a/q) &= \frac{\lambda_F\overline{\alpha_F}}{q(s-1)^2} + \frac{\lambda_F\overline{\alpha_F}(\overline{\beta_F(q/a)} + 2\log\frac{a}{q})}{q(s-1)} + \dots \\
&= \frac{\alpha_F(a/q)}{(s-1)^2} + \frac{\alpha_F(a/q)(\beta_F-2\log q)}{s-1} + \dots,
\end{split}
\]
and the result follows. \fine

\medskip
{\bf Lemma 5.4.} {\sl Let $F\in\S_2^\sharp$ with $q_F=1$. Then for $p$ prime and} $\chi$ (mod $p$) {\sl with $\chi\neq \chi_0$, $F(s,\chi)$ is holomorphic at $s=1$.}

\medskip
{\sl Proof.} From Lemmas 5.2 and 5.3 we see that the Laurent coefficients of order $-1$ and $-2$ of $F(s,a/q)$ do not depend on $a$. Hence, with obvious notation, for $\chi$ (mod $p$), $\chi\neq \chi_0$, we have 
\[
F(s,\chi) = \frac{1}{\tau(\bar{\chi})} \sum_{a=1}^p \overline{\chi(a)} F(s,-a/p) = \frac{1}{\tau(\bar{\chi})} \sum_{a=1}^p \overline{\chi(a)} \big\{\frac{A(p)}{(s-1)^2} + \frac{B(p)}{s-1} + \dots\big\},
\]
and the result follows by the orthogonality of Dirichlet characters. \fine

\bigskip
To prove Theorem 4 we choose $a=p-1$ in \eqref{2-5} of Lemma 2.2 and use Lemma 5.4 to obtain
\[
F(s,1/p) = F(s)\big(1-\frac{p}{p-1}\frac{1}{F_p(s)}\big) + G(s)
\]
with $G(s)$ holomorphic at $s=1$. Multiplying both sides of the last expression by $(s-1)^{m_F}$ and letting $s\to1$ we get
\begin{equation}
\label{5-1}
\alpha_F(1/p) = \big(1-\frac{p}{p-1}\frac{1}{F_p(1)}\big) \alpha_F.
\end{equation}
On the other hand, by Lemma 5.2 we have $\alpha_F(1/p) = \alpha_F/p$, hence comparing with \eqref{5-1} and recalling that $\alpha_F\neq0$ we obtain
\begin{equation}
\label{5-2}
F_p(1) = \big(1-\frac{1}{p}\big)^{-2}.
\end{equation}
Moreover, by the Corollary in the Introduction we have
\[
F_p(s) = \prod_{j=1}^{\partial_p}\big(1-\frac{\alpha_j(p)}{p^s}\big)^{-1}
\]
with $\partial_p\leq 2$ and $|\alpha_j(p)|\leq 1$. Therefore
\begin{equation}
\label{5-3}
|F_p(1)| = \prod_{j=1}^{\partial_p} |\sum_{m=0}^\infty \frac{\alpha_j(p)^m}{p^m}| \leq \prod_{j=1}^{\partial_p} \sum_{m=0}^\infty \frac{|\alpha_j(p)|^m}{p^m} \leq (1-\frac{1}{p})^{-\partial_p} \leq (1-\frac{1}{p})^{-2}.
\end{equation}
Comparing \eqref{5-2} and \eqref{5-3} we see that \eqref{5-2} holds if and only if $\partial_p=2$ and $\alpha_j(p)=1$ for $j=1,2$, and Theorem 4 follows.

\vskip2cm

\ifx\undefined\bysame{poly}.
\newcommand{\bysame}{\leavevmode\hbox to3em{\hrulefill}\ ,}
\fi

\vskip 1cm
\noindent
Jerzy Kaczorowski, Faculty of Mathematics and Computer Science, A.Mickiewicz University, 61-614 Pozna\'n, Poland and Institute of Mathematics of the Polish Academy of Sciences, 
00-956 Warsaw, Poland. e-mail: kjerzy@amu.edu.pl

\medskip
\noindent
Alberto Perelli, Dipartimento di Matematica, Universit\`a di Genova, via Dodecaneso 35, 16146 Genova, Italy. e-mail: perelli@dima.unige.it

\end{document}